\documentclass[%
 reprint,
 amsmath,amssymb,
 aps,
]{revtex4-2}
\pdfoutput=1

\usepackage{graphicx}
\usepackage{dcolumn}
\usepackage{bm}
\usepackage{subfigure}
\usepackage{tikz,caption} 
\usepackage{overpic}
\usepackage{xparse} 
\NewDocumentCommand{\MySubFigure}{m O{south east} m m}{ \begin{tikzpicture}[inner sep = 5pt] \node (graphic) {\includegraphics[width=\linewidth]{#1}}; \node at (graphic.#2) [anchor=#2, inner sep=-12pt] {\refstepcounter{subfigure}(\alph{subfigure})\ #3 \label{subfig:#4}}; 
\end{tikzpicture} }

\begin{document}

\preprint{APS/123-QED}

\title{Co-evolution of Vaccination Behavior and Perceived Vaccination Risk can lead to a Stag-Hunt like Game}

\author{Yuan Liu}
\affiliation{School of Sciences, Beijing University of Posts and Telecommunications, Beijing 100876, China.
}
\author{Bin Wu}
 \email{bin.wu@bupt.edu.cn}
\affiliation{School of Sciences, Beijing University of Posts and Telecommunications, Beijing 100876, China.
}

\date{\today}

\begin{abstract}
Voluntary vaccination is effective to prevent infectious diseases from spreading.
Both vaccination behavior and cognition of the vaccination risk play important roles in individual vaccination decision making.
However, it is not clear how the co-evolution of the two shapes the population-wide vaccination behavior.
We establish a coupled dynamics of epidemic, vaccination behavior and perceived vaccination risk with three different time scales. 
We assume that the increase of vaccination level inhibits the rise of perceived vaccination risk, and the increase of perceived vaccination risk inhibits the rise of vaccination level.
It is shown that the resulting vaccination behavior is similar to the stag-hunt game, provided that the basic reproductive ratio is moderate and that the epidemic dynamics evolves fast.
This is in contrast with the previous view that vaccination is a snowdrift like game.
Furthermore, we find that epidemic breaks out repeatedly and eventually leads to vaccine scares if these three dynamics evolve on a similar time scale. 
And we propose some ways to promote vaccination behavior, such as controlling side-effect bias and perceived vaccination costs.
Our work sheds light on epidemic control via vaccination by taking into account the co-evolutionary dynamics of cognition and behavior.
\end{abstract}

\maketitle

\section{Introduction}
Vaccination is regarded as one of the most effective ways to prevent the spread of infectious diseases.
It has successfully controlled many diseases, such as influenza, poliomyelitis, smallpox \cite{bonanni1999,keeling2007,may1991}, and it plays a vital role in controlling the present COVID-19 virus \cite{rawat2021,mcclung2020,chavda2021}.
Many previous models combined individual vaccination behavior with disease transmission, and fruitful results are obtained via evolutionary game theory \cite{bauch2004,bauch2013,perisic2009,wang2016}.
Typically it's assumed that vaccinated individuals pay the price for taking vaccination.
Unvaccinated individuals may be infected, or may keep healthy without paying any price \cite{bauch2004,fu2011,reluga2006,wu2011}.
However, voluntary vaccination is a long-standing social dilemma, and it is difficult to completely eradicate infectious diseases through voluntary vaccination \cite{may1991,bauch2005,bauch2003}.
This is because, vaccinated individuals will pay the cost, such as time, money, side effects and so on, in order to obtain immune ability, which is not only beneficial for the vaccinated, but also helpful for those who do not take vaccination.
When vaccination level is sufficiently high, the whole population gives rise to herd immunity \cite{quilici2015,kabir2019,peng2016}.
At this time, unvaccinated individuals are protected from infectious diseases by herd immunity without paying any cost.
This implies that unvaccinated individuals are tempting to free ride, i.e., choose not to vaccinate themselves but expect other individuals to vaccinate \cite{bauch2004,wu2011}.
Many effective methods are used to facilitate the population to escape from the vaccination dilemma.
For example, it has been explored how the update rules of strategy \cite{iwamura2018,arefin2020,huang2020}, topological structures within the population \cite{perisic2009,chang2020,wang2020,wei2019} and multi-strategy vaccination game \cite{bhattacharyya2010,ida2018,iwamura2016} improve vaccination level to create herd immunity.

Awareness is also crucial to the spread of infectious diseases in addition to human behavior \cite{funk2009,funk2010,ferguson2007,dai2021}.
For example, when cases of side-effect of vaccination are reported, the perceived vaccination risk is dramatically increasing by the communications among acquaintances \cite{bauch2004,bauch2005,reluga2006}, although the risk of the vaccination is not necessarily variant.
Whether individuals take vaccination or not is thus affected by the perceived vaccination risk \cite{galvani2007,chapman1999}. 
Studies have shown that individuals are more likely to be vaccinated, if the perceived side effects are low, or individuals think that infectious diseases are more serious \cite{brewer2007,ibuka2014}. 
Individuals who choose not to be vaccinated may be affected by the outbreak of other infectious diseases \cite{lieu2015}. 
At the same time, they will also be misled by false messages. 
For example, the statement that MMR vaccine triggered autism leads to measles outbreaks in western countries \cite{hussain2018}. 
Many individuals lack vaccine related knowledge and are vulnerable to such information. 
Even those who are in favor of vaccination may be confused by the boisterous false information, leading them to question their choice \cite{galvani2007}.
In other words, the perceived vaccination risk is an evolving opinion.
Thus, to vaccinate or not changes the social environment of the follow-up vaccination campaign, that is, perceived vaccination risk. 
In turn, the perceived vaccination risk affects individual's vaccination strategy choice in the subsequent epidemic season.

In addition to infectious disease, general ideas, sentiments and information are also `socially' infectious, in the sense that they can be driven by others' \cite{bauch2003}.
For the perception of vaccination risk, individuals are likely to take the vaccination as less secure, if few others take vaccination.
This is because popular opinions are likely to be taken as a safe one.
Previous studies do not take into account the coupling between vaccination and perceived vaccination risk, although it is known the coupling facilitates understanding of general complex systems \cite{bauch2013,chen2019}.
At present, studies have linked perceived payoffs, perceived vaccination risks and vaccination \cite{feng2018,usherwood2021} without paying attention to the feedback between the two dynamics.
How does the co-evolution of vaccination behavior and perceived vaccination risk affect the vaccination uptake level, and how to improve the vaccination uptake level if perceived risk of vaccination is also evolving?
To this end, we couple the vaccination dynamics with individual cognition of the vaccination risk, and establish a toy model to study this issue.

\section{Model}
Epidemic, vaccination behavior and perceived vaccination risk are evolving with different time scales.
Let's assume that the epidemic evolves at the unit rate,
the evolution rates of vaccination behavior and that of perceived vaccination risk relative to epidemic are $\epsilon_1$ and $\epsilon_2$, respectively.
In fact, epidemic evolves fast and the infected individuals increase exponentially \cite{Virk_2019}.
It thus implies that individuals update strategies and opinions less often, i.e., $0 < \epsilon_1 \leq 1$ and $0 < \epsilon_2 \leq 1$.

\textbf{Epidemic dynamics:}
Let's consider an infinite and well-mixed population.
We assume that the epidemic follows the classic Susceptible-Infected-Recovered (SIR model) for infectious diseases including influenzas and COVID-19.
There are three compartments in the population: susceptible (S), who can get the disease if exposed to infected individuals; infective (I), who are infected and can infect others; recovered (R), who are vaccinated or recovered from the disease and gain immunity against the disease. 
At any time, a susceptible individual gets infected, provided it interacts with an infected before it dies.
\begin{equation}
	\begin{cases}
		\dot{S}=\mu (1-x) -\beta S I -\mu S,\\
		\dot{I}=\beta S I -\gamma I -\mu I,\\
		\dot{R}=\mu x +\gamma I -\mu R,
	\end{cases}
\end{equation}
where $\mu$ is birth rate (death rate), $\beta$ is transmission rate, $\gamma$ is recovery rate, and $x$ denotes the proportion of vaccinated individuals.
In this case, the infection probability of susceptible individuals is $f(x, t) = \beta  I_t / (\beta  I_t + \mu)$.

\begin{figure}
	\includegraphics[scale=0.5]{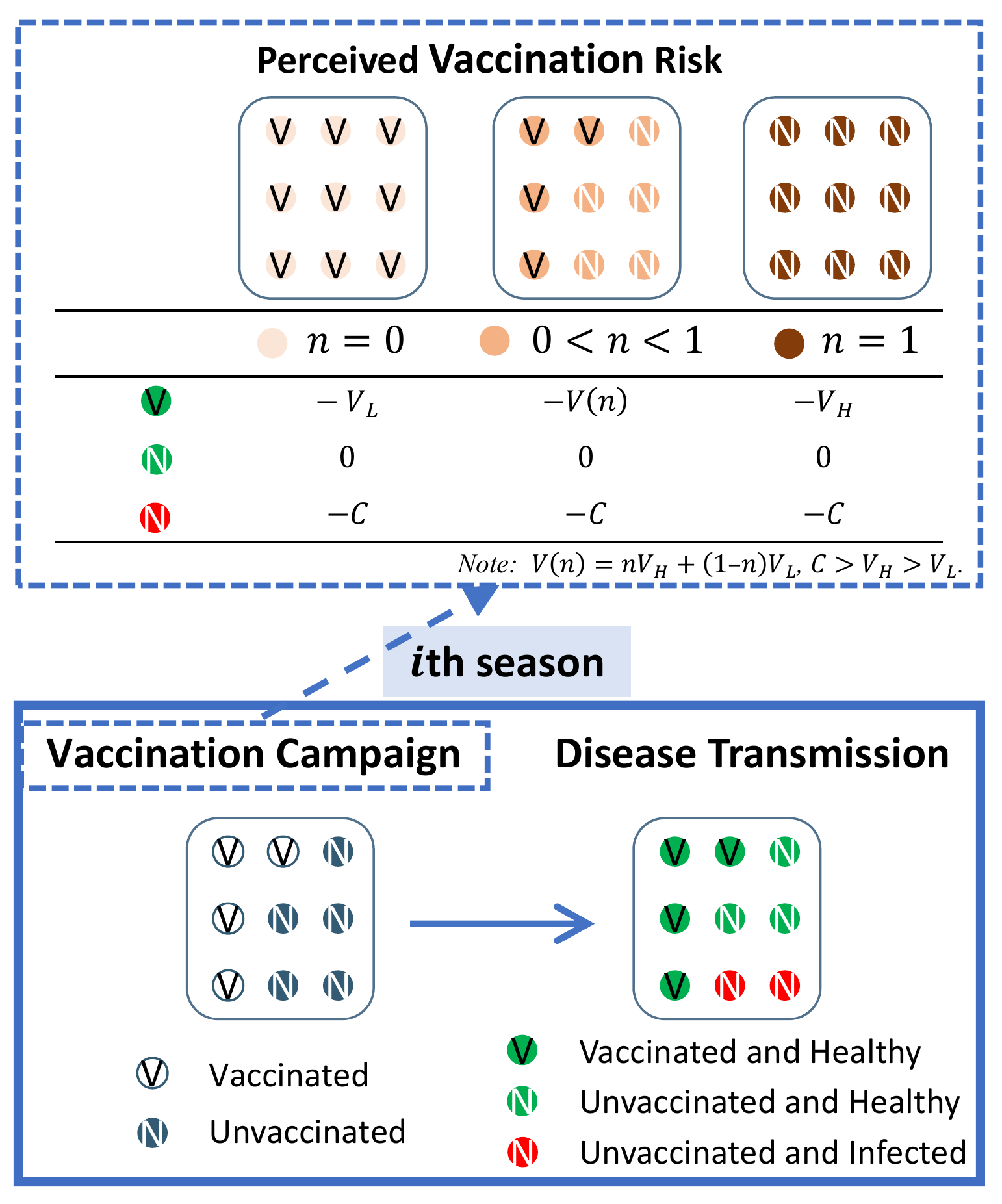}
	\caption{\label{Model}The schematic of the model. 
		The vaccination game consists of two stages: vaccination campaign and disease transmission.
		In the first stage, all the individuals have the same cognition of vaccination risk.
		They follow the majority rule in opinion dynamics: the more individuals are taking vaccination, the more secure the vaccination is regarded as. 
		At the same time, each individual has two strategies, vaccination or not.
		Vaccinated individuals will not be infected in the upcoming season (perfect vaccination).
		The unvaccinated individuals may keep healthy without paying any cost.
		However, if the unvaccinated are infected, they will pay more for recovering.
		The individual imitates other's strategy via the Fermi update rule.
		In the second stage, epidemic follows the classic SIR model. 
		There are three types of individuals at the end: vaccinated and healthy individuals, unvaccinated and healthy individuals, as well as unvaccinated and infected individuals.
		Noteworthily, the cost of vaccination is determined by the public opinion of vaccination risk.
	}
\end{figure}

\textbf{Vaccination dynamics:}
In every season of infectious diseases including influenzas and COVID-19, each individual has two strategies, vaccination or not. 
We assume that the vaccine is perfect and vaccinated individuals will not be infected in the upcoming season. 
In the vaccination stage, all the individuals have the same cognition of vaccination risk, denoted as $n$ ($0 \leq n \leq 1$).
If $n = 1$, all the individuals take the view that the vaccination risk is high.
If $n = 0$, the perceived vaccination risk is low.
In general, the vaccinated individuals pay the cost of  $V(n) = nV_H + (1-n) V_L > 0$.
Here, $V_H$ is the perceived cost of vaccination when the vaccination risk is at its maximum ($n = 1$), and $V_L$ is the perceived cost of vaccination when the vaccination risk is at its minimum ($n = 0$).
It holds that $V_H > V(n) > V_L > 0$ ($0 < n < 1$).
The unvaccinated individuals do not pay any cost.
However, if the unvaccinated are infected, they will pay the cost of $C$ for recovering, which includes the suffering from the disease as well as social economic costs.
Typically, the cost of recovery is larger than the highest perceived vaccination cost, that is, $C > V_H$ (See Fig.~\ref{Model}).

In the stage of disease transmission, the vaccinated individuals are not at risk of infection.
The infection probability of unvaccinated individuals is $f(x,t)$, in which $x$ is the proportion of vaccinated individuals.
There are three types of individuals in the population.
Vaccinated and healthy individuals' proportion is $x$ with payoff $- V (n)$;
Unvaccinated and infected individuals' proportion is $(1-x) f(x,t)$ with payoff $- C$;
Unvaccinated and healthy individuals' proportion is $(1-x) (1-f(x,t))$ with payoff $0$.
At this time, the individual imitates the strategy of other individuals via the Fermi update rule \cite{fu2011, traulsen2005, wu2010}, that is, $i$ and $j$ are randomly selected, and the probability of $i$ learning $j$' strategy is $[1+\exp{(-\beta (\pi_j-\pi_i))}]^{-1}$.
Here, $\pi_i$ and $\pi_j$ are the perceived payoffs of $i$ and $j$ respectively, $\beta \ge 0$ is selection intensity, which determines how much the payoff difference affects the individual's strategy.
When the selection intensity is zero, it implies that individual $i$ learns $j$'s strategy with probability of one half, regardless of the payoffs of the two individuals.
When the selection intensity is strong, i.e., $\beta\to +\infty$, it implies that individual $i$ is almost sure to learn $j$'s strategy, even though the payoff of $j$'s strategy is only a slightly higher than that of $i$'s strategy.
Noteworthily, the selection intensity mirrors the inverse temperature in statistical physics \cite{TRAULSEN2007522, ZHANG2021126543, traulsen2006}. 
Thus, vaccination dynamics is described by:
\begin{eqnarray}
	\dot{x}=&&x(1-x)[(1-f(x,t))\tanh{(\frac{\beta}{2} (-V(n)))}\\
  &&+f(x,t)\tanh{(\frac{\beta}{2}(-V(n)+C))}].
\end{eqnarray}

\textbf{Perceived vaccination risk dynamics:}
In every epidemic season,
the vaccination level of the current season will have an impact on the individual's perceived vaccination risk during the vaccination campaign.
We follow the majority rule in opinion dynamics
in the sense that
the more individuals are taking vaccination, the more secure the vaccination is regarded as.
This is because the social groups are likely to adopt popular opinions,
and the popular opinions are taken as less harmful.
As a basic toy model, we propose that
i) the risk is between $0$ and $1$;
ii) the more vaccinators there are, the less the risk of the vaccination is taken;
iii) the more unvaccinators there are, the more the perceived risk of the vaccination is.
To this end, we assume that the greater the side-effect bias $\theta$ is, the stronger role the unvaccinated play in the dynamics of the perceived vaccination risk.
\begin{equation}
	\dot{n}=n(1-n)[-x+(1+\theta) (1-x)],
\end{equation}
where $\theta>0$ is side-effect bias, which measures the ratio of the enhancement rates to degradation rates of vaccinated and unvaccinated individuals, respectively.
The logistic term $n(1-n)$ ensures that the perceived vaccination risk is restrained to $[0,1]$, as it is defined within. 
In addition, the term $(-x+(1+\theta) (1-x))$ describes the side-effect bias with the assumption that the more vaccinators there are, the less the perceived risk of the vaccination is; the more unvaccinators there are, the more the perceived risk of the vaccination is  (see Appendix \ref{dynamic equations} for details).
In fact, the side-effect bias could be controlled by some measures.
For example, side-effect bias decreases by reporting the positive utility and data of vaccinated individuals \cite{Prasad2022, Plumb2022}. 
Otherwise, the spreading of negative, exaggerated and untrue statements increase the side-effect bias \cite{hussain2018, sanyaolu2019}.

Therefore, epidemic, vaccination behavior and perceived vaccination risk are evolving with time in our model. 
However, the evolution rates of these three dynamics may be different.
Under weak selection limit $\beta \to 0^+$, after a time rescaling which does not change the asymptotic dynamics, then the resulting coupled dynamics of epidemic described by SIR model, vaccination behavior and perceived vaccination risk are given by: 
\begin{eqnarray}
	\begin{cases}
		\dot{S}=\mu (1-x) -\beta S I -\mu S,\\
		\dot{I}=\beta S I -\gamma I -\mu I,\\
		\dot{R}=\mu x +\gamma I -\mu R,\\
		\dot{x}=\epsilon_1 x(1-x)[(f(x,t) C- V(n)],\\
		\dot{n}=\epsilon_2 n(1-n)[-x+(1+\theta) (1-x)],
	\end{cases}
	\label{eq2}
\end{eqnarray}
In spite of its simplicity, our model captures i) the change of vaccination level affects the vaccination risk; ii) the individual's cognition of vaccination risk affects whether to vaccinate or not in the next epidemic season.

\section{Dynamical Analysis}
In this section, we consider two cases:
i) the epidemic evolves much faster than vaccination behavior and perceived risk;
ii) these three dynamics evolve with a similar time scale.
We find that co-evolutionary vaccination is similar to a stag-hunt game in case i), and repeated outbreaks in epidemic leads to vaccine scare in case ii).

\subsection{Co-evolution of Vaccination Behavior and Perceived Vaccination Risk can lead to a Stag-Hunt like Game}
There are some studies on coevolutionary dynamics of epidemics and vaccination \cite{bhattacharyya2010,cui2017,PhysRevE.88.032803,rifhat2021}.
Repeated outbreaks of epidemic are present due to external factors, such as immigration, zoonoetic and thing-to-human, even if the epidemic dynamics has reached equilibrium \cite{madhav2017,cohen2005}.
Then a new framework is proposed \cite{fu2011,alessio2013,bauch2012}, in which the epidemic transmission and vaccination campaign are divided into two stages.
Once the epidemic dynamics reached its stationary regime, individuals choose whether or not to vaccinate before the next season begins. 
It resembles a time-scale separation, which is worth considering.
Based on similar considerations, we assume the epidemic evolves much faster than vaccination behavior and perceived risk.
It implies that once the epidemic dynamics dies out, that is, at the final state of epidemic season, individuals choose whether or not to vaccinate before the next outbreak begins.
At the same time, we assume that vaccination behavior evolves at the same rate as perceived vaccination risk.
Thus, we set $\epsilon_1, \epsilon_2 \ll 1$ to indicate that the epidemic evolves much faster than vaccination behavior and perceived risk.
For simplicity, we assume $\epsilon_1 = \epsilon_2$.
In this case, the infection probability of susceptible individual $f(x)$ is seen as the limit of $f (x, t) = \beta  I_t / (\beta    I_t + \mu )$ when $t$ tends to infinity, which shows that the epidemic dynamics has reached the equilibrium. 
As mentioned before, infection probability $f(x)$ is given by \cite{wu2011,traulsen2006,bauch2005}:
\begin{equation}
	f(x)=\begin{cases}
		1-\frac{1}{R_0 (1-x)}, & \mbox{if } 0 \leq x < 1-\frac{1}{R_0} \\
		0, & \mbox{if } x \ge 1-\frac{1}{R_0}
	\end{cases},
\end{equation}
if the epidemic follows the classic SIR model.
Here, $R_0$ is the basic reproduction ratio, which is the average infection number of secondary infections among infected individuals in the susceptible individuals.
Under weak selection limit $\beta \to 0^+$, then the coupled dynamics is given by:
\begin{equation}
	\begin{cases}
		\dot{x}=x(1-x)[f(x)C-V(n)],\\
		\dot{n}=n(1-n)[-x+(1+\theta) (1-x)].
	\end{cases}
	\label{eq1}
\end{equation}

The above equation is much less complicated than Eq.~\ref{eq2}, which can facilitate us to provide analytical insight.

We find that there are only two stable regimes when the basic reproductive ratio $R_0$ is moderate.
One stable regime has a high vaccination level with the lowest perceived vaccination risk.
The other stable regime has low vaccination level with highest perceived vaccination risk.
We show how side-effect bias, disease severity and individual payoff alter vaccination level and cognition risk of vaccination to shed light on epidemic control.

\subsubsection{Two stable regimes}
When $R_0$ is too small or too large, there is only one stable fixed point.
For small basic reproductive ratio, $R_0 < 1$, the disease disappears naturally and the infection risk ends up is low.
Individuals believe that the vaccination risk is higher than that of infection.
Rational individuals do not choose to vaccinate in the following vaccination campaign.
Based on our model, the perceived vaccination risk increases with the decrease of vaccination level.
Thus, the vaccination level reaches zero and the perceived vaccination risk reaches one.
When $R_0\gg1$, however, the number of infected individuals increases exponentially and the disease lasts for a long time.
At this time, individuals believe that the infection risk is very high.
In other words, the cost and side effects caused by vaccination are obviously not worth mentioning.
In this case, the vaccination level maintains at a level $1-C/[R_0(C-V_L)]$, and the vaccination risk reaches zero (See Appendix \ref{fixed points} for details).
Noteworthily, $1-C/[R_0(C-V_L)]$ is the vaccination level at which the perceived vaccination risk is not evolving and keeps the lowest.

However, the dynamics is completely different when the basic reproductive ratio is moderate.
When
\begin{eqnarray*}
 	\frac {(2+ \theta) C} {C-V_L} < R_0 < \frac {(2+ \theta) C} {C-V_H},
\end{eqnarray*}
two and only two stable regimes are present.
And if
\begin{eqnarray*}
	  \frac{(2+\theta)C}{C-V_L}>\frac{C}{C-V_H},
\end{eqnarray*}
the two stable fixed points are
\begin{eqnarray*}
 	 (1)(1-\frac{C}{R_0(C-V_L)},0),
\end{eqnarray*}
\begin{eqnarray*}
	 (2)(1-\frac{C}{R_0(C-V_H)},1).
\end{eqnarray*}
The phase diagram satisfying these conditions is shown in Fig.~\ref{PD1}. %
\begin{figure}
	\includegraphics[scale=1]{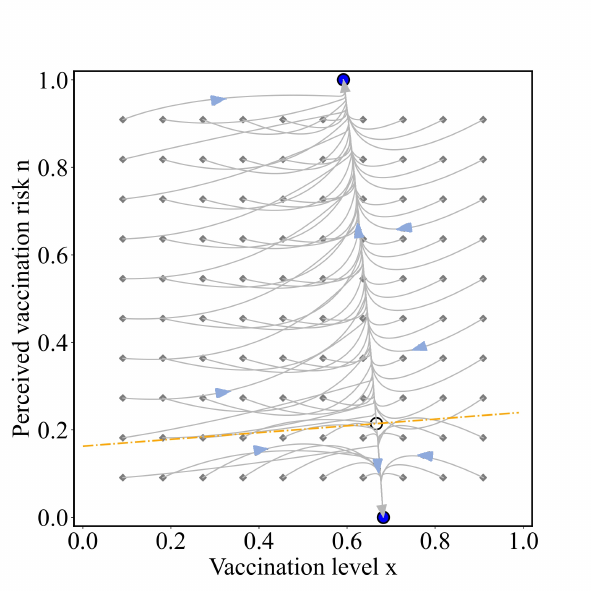}
	\caption{\label{PD1}Phase plane dynamics of $x-n$. 
	The gray curves denote the realized orbits and the arrows denote the direction of the trajectories. 
	The black open circle denotes the unstable internal fixed point. 
	The blue solid circles denote the stable point which has high vaccination level with the lowest perceived vaccination risk, and the other stable point which has low vaccination level with highest perceived vaccination risk, respectively. 
	The orange dotted line separates the attraction domain of the two stable fixed points  approximately.
	Parameters: $R_0=3.5, C=10, V_H=3, V_L=1, \theta=1$. 
}
\end{figure}
The two stable fixed points refer to high vaccination level ($1-C/[R_0(C-V_L)]$) with the lowest perceived vaccination risk and low vaccination level ($1-C/[R_0(C-V_H)]$) with the highest perceived vaccination risk, respectively.
This dynamical outcome is absent if the perceived vaccination risk is not co-evolving with behavior \cite{bauch2004,bauch2005,fu2011,wu2011,reluga2006}.
Thus, we are more interested in the case that the basic reproductive ratio is moderate (see Appendix \ref{fixed points}).

From the perspective of perceived vaccination risk, on one hand, when the initial perceived vaccination risk is high, the vaccination level will be reduced.
Lower vaccination levels, in turn, give rise to higher perceived vaccination risks.
When the time is sufficiently long, the perceived vaccination risk reaches one, but the vaccination level remains at a low level $1-C/[R_0(C-V_H)]$.
On the other hand, when the perceived vaccination risk in the initial state is low, the vaccination level increases.
The increase of vaccination level leads to lower perceived vaccination risk.
Eventually, the perceived vaccination risk reaches zero, and the vaccination level remains at a high level $1-C/[R_0(C-V_L)]$.

From the perspective of vaccination level, if the initial vaccination level is high, rational individuals choose not to vaccinate.
It reduces the vaccination level and increase the perceived vaccination risk.
Finally, the perceived vaccination risk reaches one, and the vaccination level remains at a low level $1-C/[R_0(C-V_H)]$.
If the initial vaccination level is low, individuals think that they have a high risk of infection, as a result, they choose to vaccinate.
The increase of vaccination level reduces the perceived vaccination risk.
Then the perceived vaccination risk reaches zero, and the vaccination level remains at a high level $1-C/[R_0(C-V_L)]$.
Therefore, when $R_0$ is moderate, there are two stable fixed points, one of which has a higher vaccination level and the lowest vaccination risk and the other one is versa.

\subsubsection{Mechanisms to promote vaccination behavior: a stag-hunt game perspective}
Since the vaccination level of fixed point (1) is higher and the perceived vaccination risk is the lowest, we are more concerned about how to enlarge its attraction basin.
We regard the vaccination level to be improved, if the attraction basin of fixed point (1) is larger.
To obtain higher vaccination level, we discuss how side-effect bias, vaccination cost and disease severity alter vaccination level and perceived vaccination risk (see Fig.~\ref{area}).%
\begin{figure*}
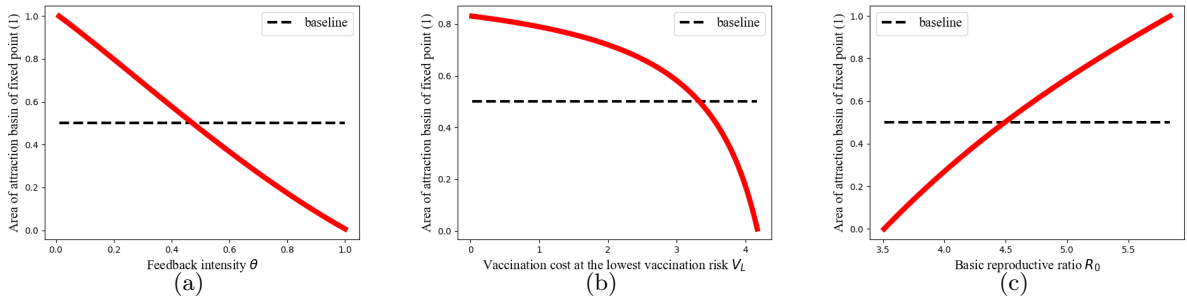
 
	\stepcounter{figure} 
	\begin{minipage}{0.3\linewidth} 
	\MySubFigure{STheta.pdf}[south]{}{STheta} 
	\end{minipage}
	\begin{minipage}{0.3\linewidth} 
	\MySubFigure{SVL.pdf}[south]{}{SVL} 
	\end{minipage}
    \begin{minipage}{0.3\linewidth}
    \MySubFigure{SR0.pdf}[south]{}{SR0}
    \end{minipage}
	\addtocounter{figure}{-1} 
	\caption{The red solid line denotes the change of the attraction basin area of fixed point (1) with side-effect bias $\theta$ (a), vaccination cost at the lowest perceived vaccination risk $V_L$ (b) and basic reproductive ratio $R_0$ (c). 
	The black dotted line denotes the baseline, and the red solid line above it showes that the regime with high vaccination level and low perceived vaccination risk is dominant. 
	Parameters are $R_0=4, C=4, V_H=2, V_L=1$ for graph (a), $R_0=3.6, C=10, V_H=5, \theta=0.1$ for graph (b) and $C=3.5, V_H=2, V_L=1, \theta=0.5$ for graph (c).}
	\label{area} 
\end{figure*} 

We find that when $ \frac {(2+ \theta) C} {C-V_L} < R_0 < \frac {(2+ \theta) C} {C-V_H}$, the dynamics is similar to that of the stag-hunt game.
The stag-hunt game has two stable Nash equilibria.
All rational individuals choose either stag or hare.
And the two stable fixed points are separated by an unstable fixed point.
Technically, the unstable fixed point captures the attraction basin of the two stable equilibria.
The equilibrium which has larger attraction basin is regarded as favored by natural selection (risk dominant strategy).
In the one-dimensional model, the internal equilibrium is sufficient to capture which strategy is dominant, that is, which strategy has a larger attraction basin.
However, in our two-dimension model, what we need to find is no longer a point, but a curve composed of countless points.
Thus, it is challenging to clarify which stable regime has a larger attraction basin.
Here, we are to give an estimation of the attraction basin of the more effective vaccination regimes, that is, low perceived risk and high level of vaccination  (see the orange dotted line in Fig.~\ref{PD1}).
We find that such a curve passes through the internal saddle point
\begin{eqnarray*}
  (3)(\frac{\theta+1}{\theta+2},-\frac{C(2-R_0+\theta)+R_0V_L}{R_0(V_H-V_L)}),
\end{eqnarray*}
in which $ \frac {(2+ \theta) C} {C-V_L} < R_0 < \frac {(2+ \theta) C} {C-V_H}$, and the direction is a stable direction of the point.
If the position of the initial point of the system is below this line, the vaccination level will reach $1-C/[R_0(C-V_L)]$ and the perceived vaccination risk will reach zero eventually.
If the position of the initial point is above the line, the vaccination level will reach $1-C/[R_0(C-V_H)]$ and the perceived vaccination risk will reach one eventually.
To obtain this approximated line, we calculate the negative eigenvalue and its eigenvector of the Jacobian matrix at the saddle point, then approximately divide the attraction basin of the two stable regimes through a line passing through the point and the direction is the eigenvector corresponding to the negative eigenvalue (see Appendix \ref{attraction basin}).

\textbf{Low side-effect bias promotes vaccination.}
The two equilibria has nothing to do with the side-effect bias.
But when $C$, $V_ H$, $V_ L$ and $R_0$ do not change, the attraction basin of fixed point (1) shrinks with the increase of $\theta$, as shown in the Fig.~\ref{subfig:STheta}. 
We select three representative points, namely $\theta = 0.05,~ 0.5,~ 0.8$, and give the phase diagrams when other parameters remain unchanged (see Fig.~\ref{Theta}).%
\begin{figure*}
	\stepcounter{figure} 
	\begin{minipage}{0.3\linewidth} 
	\MySubFigure{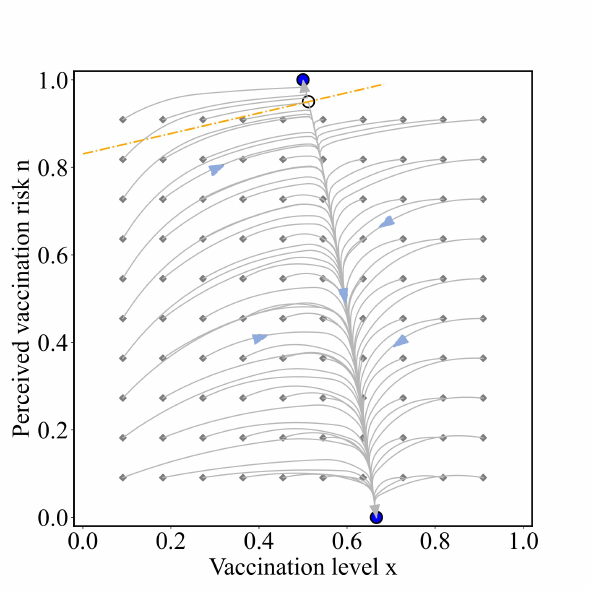}[south]{}{theta1} 
	\end{minipage}
	\begin{minipage}{0.3\linewidth} 
	\MySubFigure{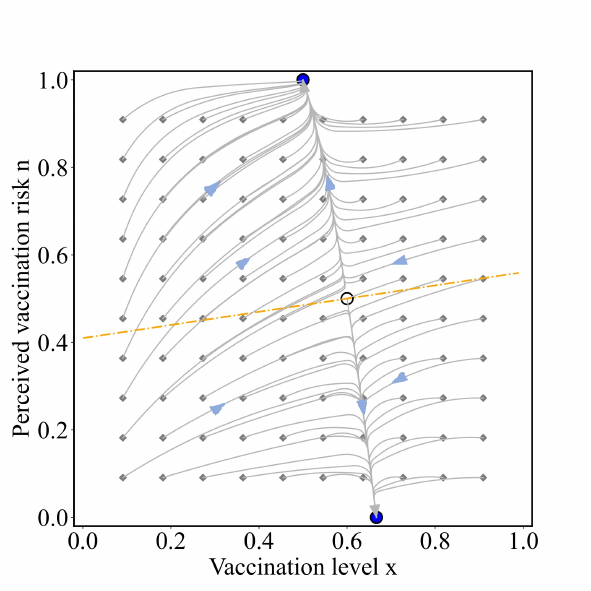}[south]{}{theta2} 
	\end{minipage}
    \begin{minipage}{0.3\linewidth}
    \MySubFigure{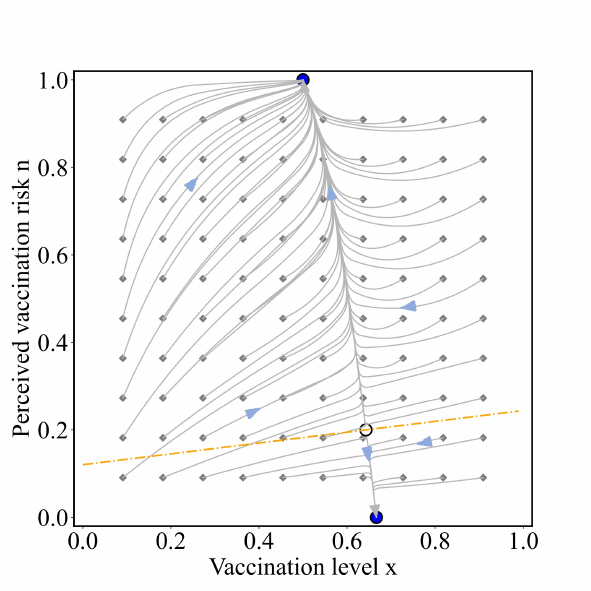}[south]{}{theta3}
    \end{minipage}
	\addtocounter{figure}{-1} 
	\caption{Phase plane dynamics of $x-n$ when $\theta=0.05 \text{ (a)},~ 0.5 \text{ (b)},~ 0.8 \text{ (c)}$.
    The gray curves denote the realized orbits and the arrows denote the direction of the trajectories. 
	The black open circle denotes the unstable internal fixed point. 
	The blue solid circles denote the stable point which has high vaccination level with the lowest perceived vaccination risk, and the other stable point which has low vaccination level with highest perceived vaccination risk, respectively. 
	The orange dotted line separates the attraction domain of the two stable fixed points  approximately.
    The areas of attraction basin of fixed point (1) (the bottom one) are about 0.9489, 0.4608 and 0.1655, respectively. 
	The parameters are the same as Fig.~\ref{subfig:STheta}.}
	\label{Theta} 
\end{figure*} 

Intuitively, when the side-effect bias increases, individuals believe that the current perceived vaccination risk is even higher once they contact an unvaccinated individual.
They tend not to take vaccination, the vaccination level thus reduces.
When the vaccination level decreases, the probability of individuals interacting unvaccinated individuals will increase.
And the perceived vaccination risk will further increase.

On the other hand, we find that the vaccination level of internal unstable fixed point (3) is $x^* = \frac {\theta + 1} {\theta + 2}$.
When $\theta > 0$, $x^*$ is a monotonically increasing of the side-effect bias $\theta$ and is always between $(0.5,1)$.
It mirrors the unstable equilibrium in stag-hunt game determines the attraction basin of the stag behavior, which is socially optimum. Similarly, here $x^{*}$ also sheds light on the attraction basin.

\textbf{Low perceived vaccination cost promotes vaccination.}
The attraction basin of fixed point (1) shrinks with the increase of vaccination cost $V_L$, as shown in Fig.~\ref{subfig:SVL}. 
We also select three representative points, namely, $V_L = 1,~ 3,~4$, and show the phase diagrams when other parameters remain unchanged (see Fig.~\ref{Vl}).%
\begin{figure*}
	\stepcounter{figure} 
	\begin{minipage}{0.3\linewidth} 
	\MySubFigure{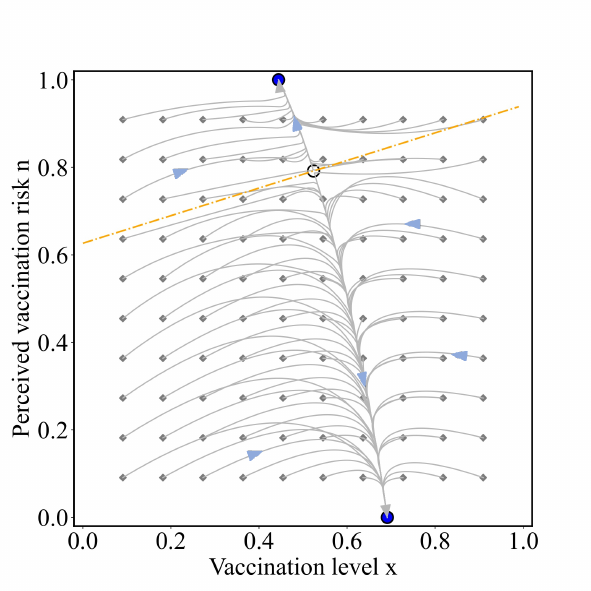}[south]{}{VL1} 
	\end{minipage}
	\begin{minipage}{0.3\linewidth} 
	\MySubFigure{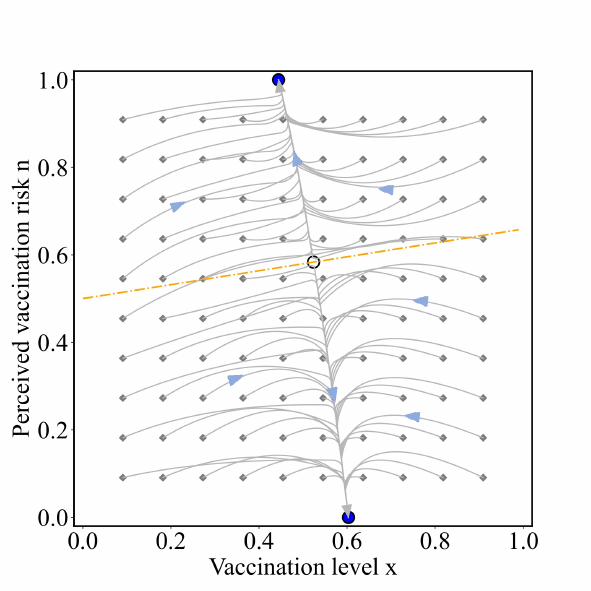}[south]{}{VL3} 
	\end{minipage}
    \begin{minipage}{0.3\linewidth}
    \MySubFigure{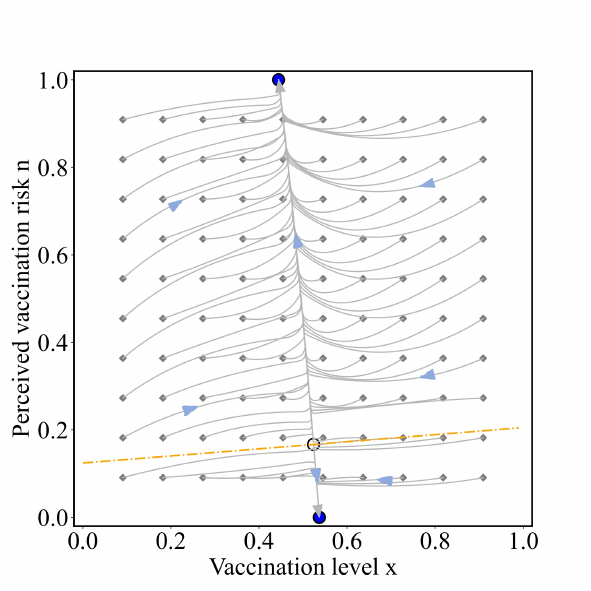}[south]{}{VL4}
    \end{minipage}
	\addtocounter{figure}{-1} 
	\caption{Phase plane dynamics of $x-n$ when $V_L=$ 1 (a), 3 (b) and 4 (c).
	The gray curves denote the realized orbits and the arrows denote the direction of the trajectories. 
	The black open circle denotes the unstable internal fixed point. 
	The blue solid circles denote the stable point which has high vaccination level with the lowest perceived vaccination risk, and the other stable point which has low vaccination level with highest perceived vaccination risk, respectively. 
	The orange dotted line separates the attraction domain of the two stable fixed points  approximately.
	The areas of attraction basin of fixed point (1) (the blue one) are about 0.7891, 0.5795 and 0.1644, respectively. 
	The parameters are the same as Fig.~\ref{subfig:SVL}.}
	\label{Vl} 
\end{figure*} 

It's shown that the whole vaccination cost $V(n)$ increases when the vaccination cost at the lowest perceived risk $V_L$ increases.
Intuitively, when the epidemic is not serious (in this case, $R_0 = 3.6$), individuals believe that the risk of vaccination outperforms the risk of disease.

\textbf{High basic reproductive ratio promotes vaccination.}
We find that when $C$, $V_ H$, $V_L$ and $\theta$ do not change, the attraction basin of fixed point (1) spans with the increase of basic reproductive ratio $R_0$, as shown in the Fig.~\ref{subfig:SR0}.
We show three representative points, namely, $R_0 = 3.6,~ 4,~ 5.5$, and give the phase diagrams respectively when other parameters remain unchanged (see Fig.~\ref{R0}). %
\begin{figure*}
	\stepcounter{figure} 
	\begin{minipage}{0.3\linewidth} 
	\MySubFigure{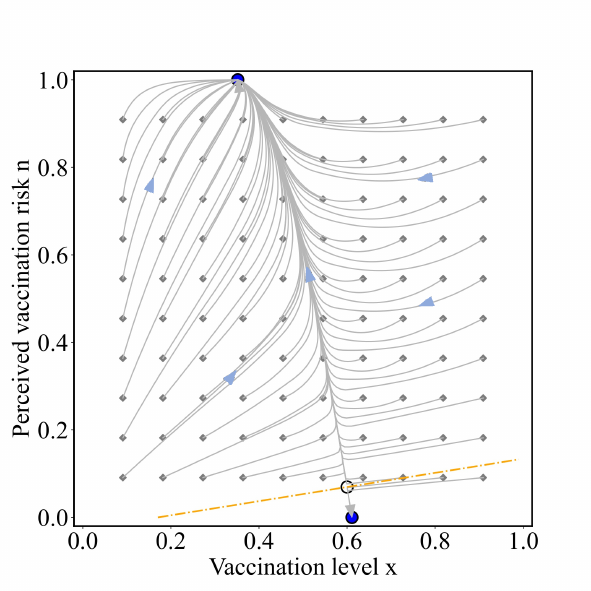}[south]{}{R01} 
	\end{minipage}
	\begin{minipage}{0.3\linewidth} 
	\MySubFigure{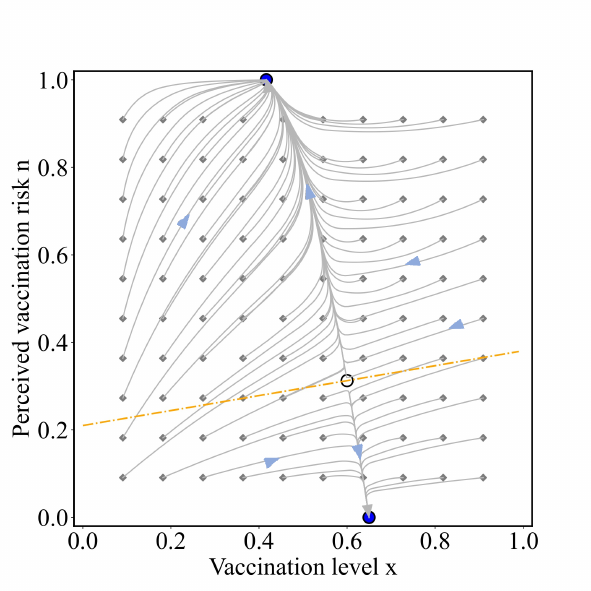}[south]{}{R02} 
	\end{minipage}
    \begin{minipage}{0.3\linewidth}
    \MySubFigure{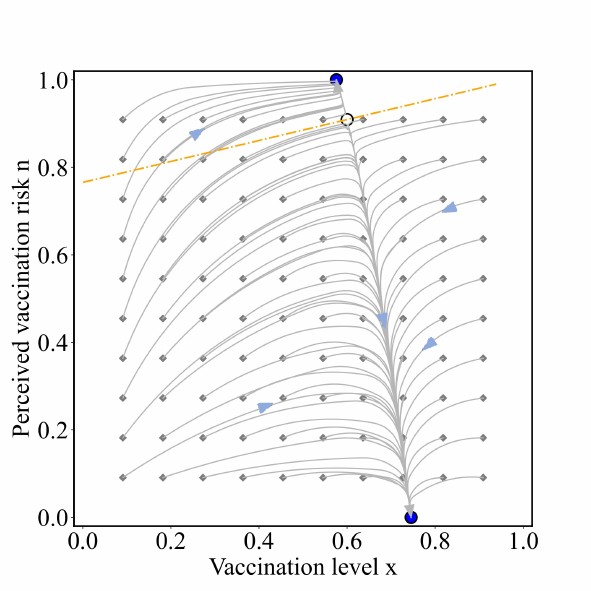}[south]{}{R03}
    \end{minipage}
	\addtocounter{figure}{-1} 
	\caption{Phase plane dynamics of $x-n$ when $R_0=$ 3.6 (a), 4 (b) and 5.5 (c).
	The gray curves denote the realized orbits and the arrows denote the direction of the trajectories. 
	The black open circle denotes the unstable internal fixed point. 
	The blue solid circles denote the stable point which has high vaccination level with the lowest perceived vaccination risk, and the other stable point which has low vaccination level with highest perceived vaccination risk, respectively. 
	The orange dotted line separates the attraction domain of the two stable fixed points  approximately.
	The areas of attraction basin of fixed point (1) (the blue one) are about 0.0586, 0.2743 and 0.8885, respectively. 
	The parameters are the same as Fig.~\ref{subfig:SR0}.}
    \label{R0} 
\end{figure*}

If the measures of disease prevention (including time cost, cognition of disease risk, etc.) do not change, the basic reproductive ratio $R_0$ is large, which implies the epidemic is serious. 
In this case, individuals believe that the perceived vaccination risk is significantly reduced compared with the infection risk.
Thus, rational individuals tend to take vaccination.
According to our assumption, the vaccination level increases then the perceived vaccination risk decreases.
With the increase of $R_0$, the area of attraction basin of the fixed point (1) increases.
The whole vaccination level is high and perceived vaccination risk is the lowest. 
It implies that more and more individuals take the view that the perceived vaccination risk is low and then choose to vaccinate.

\subsection{Repeated outbreaks in epidemic can lead to vaccine scare}
In fact, disease transmission and vaccination campaigns are difficult to delineate precisely into two separate stages.
For example, the number of COVID-19 infections never fully cleared in the past two years.
It implies that the epidemic has not yet reached its equilibrium.
During the disease transmission stage, individuals change their vaccination behavior and perceived vaccination risk.
Thus, we consider the model by incorporating three time scales for the three dynamics we take into account (SIR, vaccination behavior and risk assessment), i.e. Eq.~\ref{eq2}.
We aim to investigate how the relative time scales of epidemics and vaccination alters the outcome. 
At any time $t$, the infection probability of unvaccinated susceptible individuals is $f(x, t) = \beta  I_t / (\beta  I_t + \mu)$, in which $x$ denotes the proportion of vaccinated individuals. 
Here, we assume that initial proportion of infected individuals is $I_0=0.1$, initial recovered individuals is zero and initial proportion of susceptible individuals is $S_0 = 1 - I_0 - x_0$, in which $x_0$ denotes the initial proportion of vaccinated individuals.

\textbf{If epidemic evolves much faster, the time scale separation facilitates us to analyze the coupled dynamics.}
For Eq.~\ref{eq2}, the basic reproduction ratio is $R_0 = \beta / (\gamma + \mu)$. 
When the basic reproductive ratio is moderate, that is, $R_0$ is in $(\frac {(2+ \theta) C} {C-V_L},\frac {(2+ \theta) C} {C-V_H})$ approximately,
two and only two stable regimes are present.
One stable regime has high vaccination level and low perceived vaccination risk and
the other one has low vaccination level and high perceived vaccination risk (See Fig.~\ref{SIRa}). 
\begin{figure*}
	\stepcounter{figure} 
	\begin{minipage}{0.45\linewidth} 
	\MySubFigure{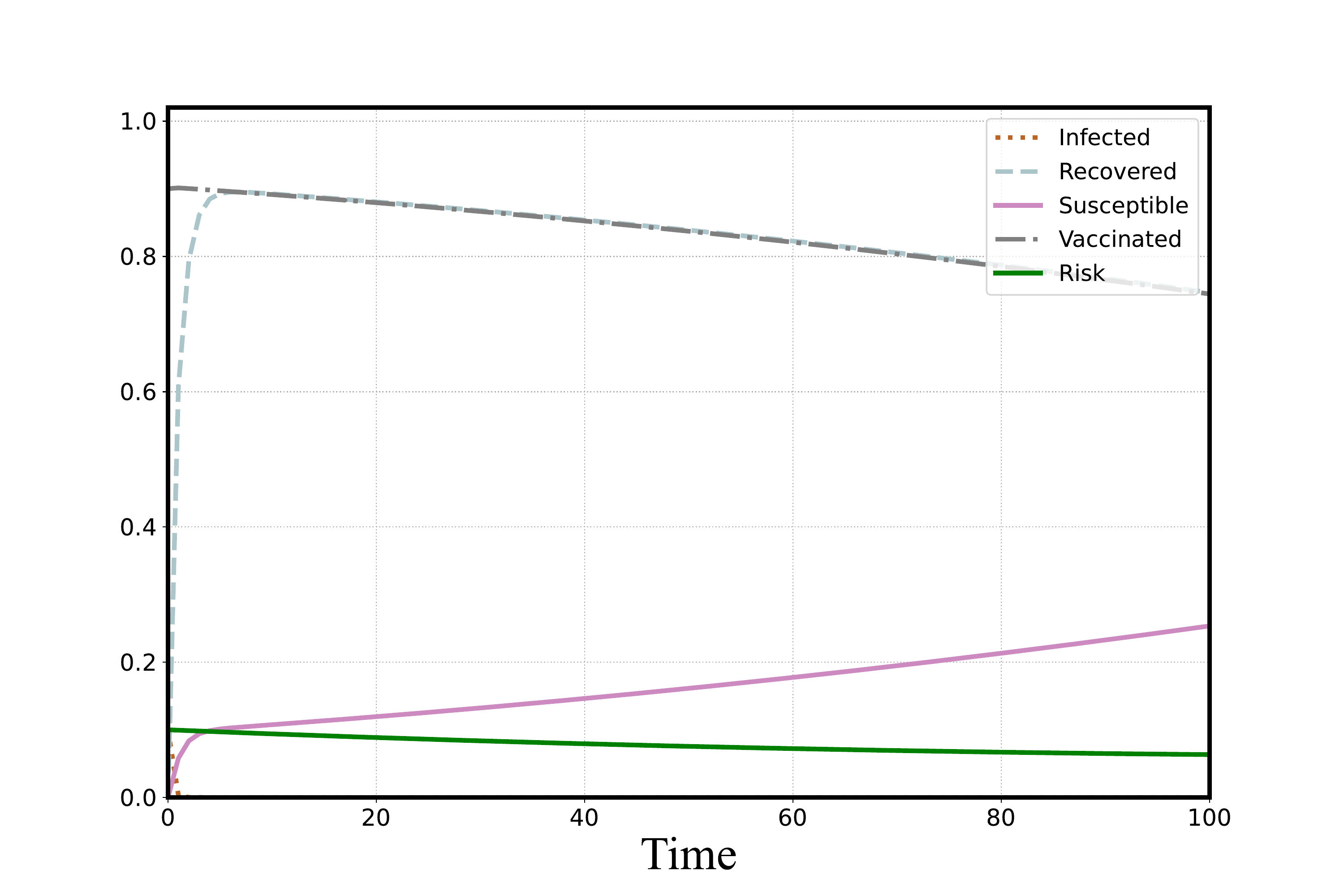}[south]{$x_0 = 0.9, n_0 = 0.1$}{SIRa-1} 
	\end{minipage}
	\begin{minipage}{0.45\linewidth} 
	\MySubFigure{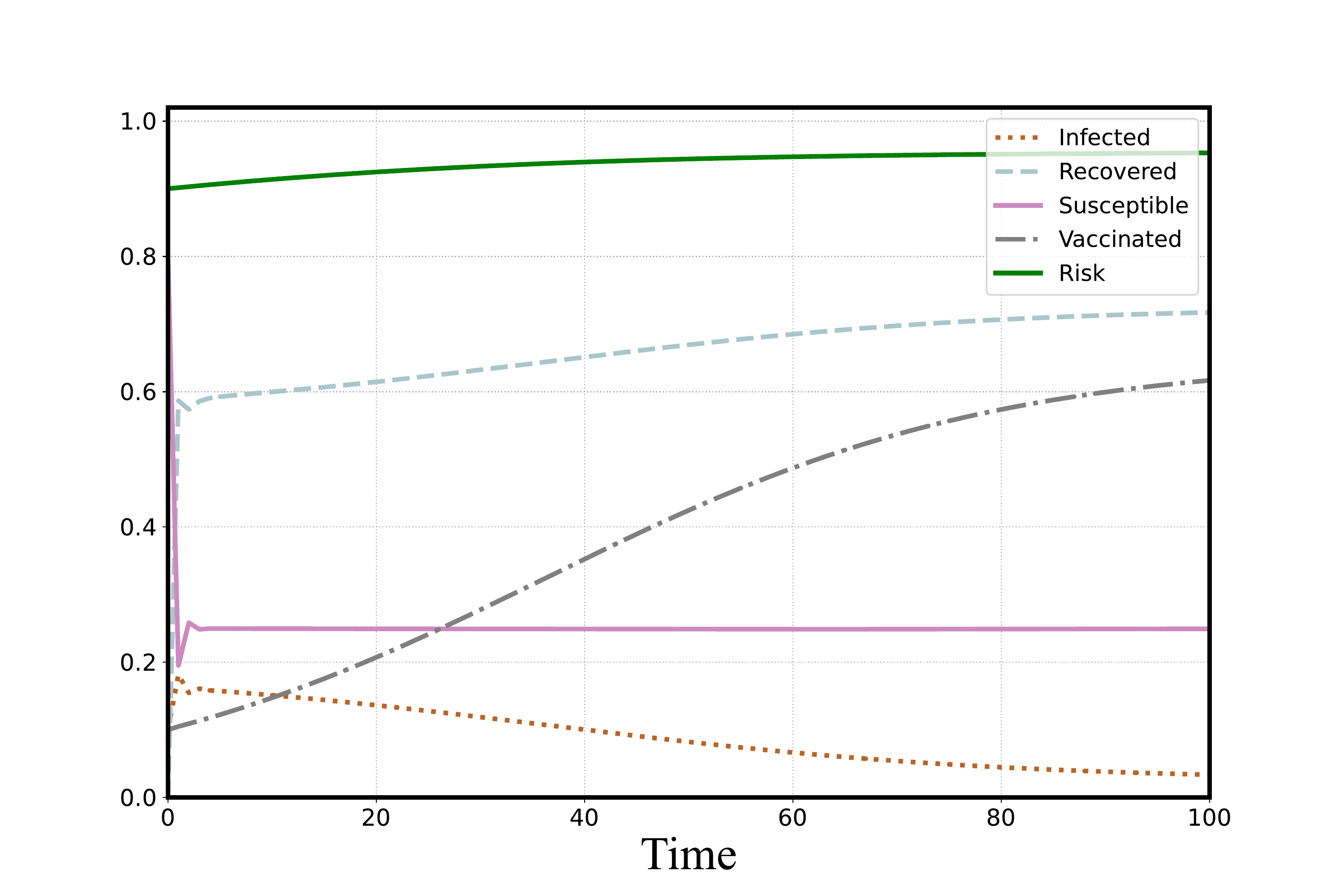}[south]{$x_0 = 0.1, n_0 = 0.9$}{SIRa-2} 
	\end{minipage}
	\addtocounter{figure}{-1} 
	\caption{\label{SIRa} There are two stable regimes of vaccination, and the initial state influences the end of coupled dynamics. (a) Initial vaccination level $x_0 = 0.9$ and perceived vaccination risk $n_0=0.1$, and the final vaccination level reach a higher proportion and perceived risk reach zero. 
	(b) Initial vaccination level $x_0=0.1$ and perceived vaccination risk $n_0 = 0.9$, and the final vaccination level reach a lower proportion and perceived risk reach one. 
	The result is consistent with the case given by Eq.~\ref{eq1}. 
	Parameters: $\epsilon_1 = \epsilon_2 = 0.01, \mu=1, \beta=16, \gamma=3, \theta=1, C=10, V_H=3, V_L=1$.
}
\end{figure*}
It is consistent with the analysis.
If both $\epsilon_1$ and $\epsilon_2$ are much smaller than 1, it degenerates to the case given by Eq.~\ref{eq1}, in which time scale separation is applied to obtained analytical insights.

\textbf{If these three dynamics evolve at the similar time scale, it can lead to vaccine scare.}
In this case, the epidemic evolves with vaccination behavior and perceived risk at the similar time scale.
Once epidemic outbreaks, infected individuals increase. 
Then parts of susceptible individuals choose to vaccinate, infected individuals gradually decrease in number. 
However, epidemic is not serious, it does not die out naturally and the vaccination level does not lead to herd immunity.
The proportion of infected individuals does not reach zero, and the epidemic continues.
Thus, outbreaks of epidemic come and go, and the susceptible, the infected, the recovered and the vaccinated individuals fluctuate in number. 
On the other hand, the perceived vaccination risk reaches one (See Fig.~\ref{SIRb}).
\begin{figure*}
	\stepcounter{figure} 
	\begin{minipage}{0.45\linewidth} 
	\MySubFigure{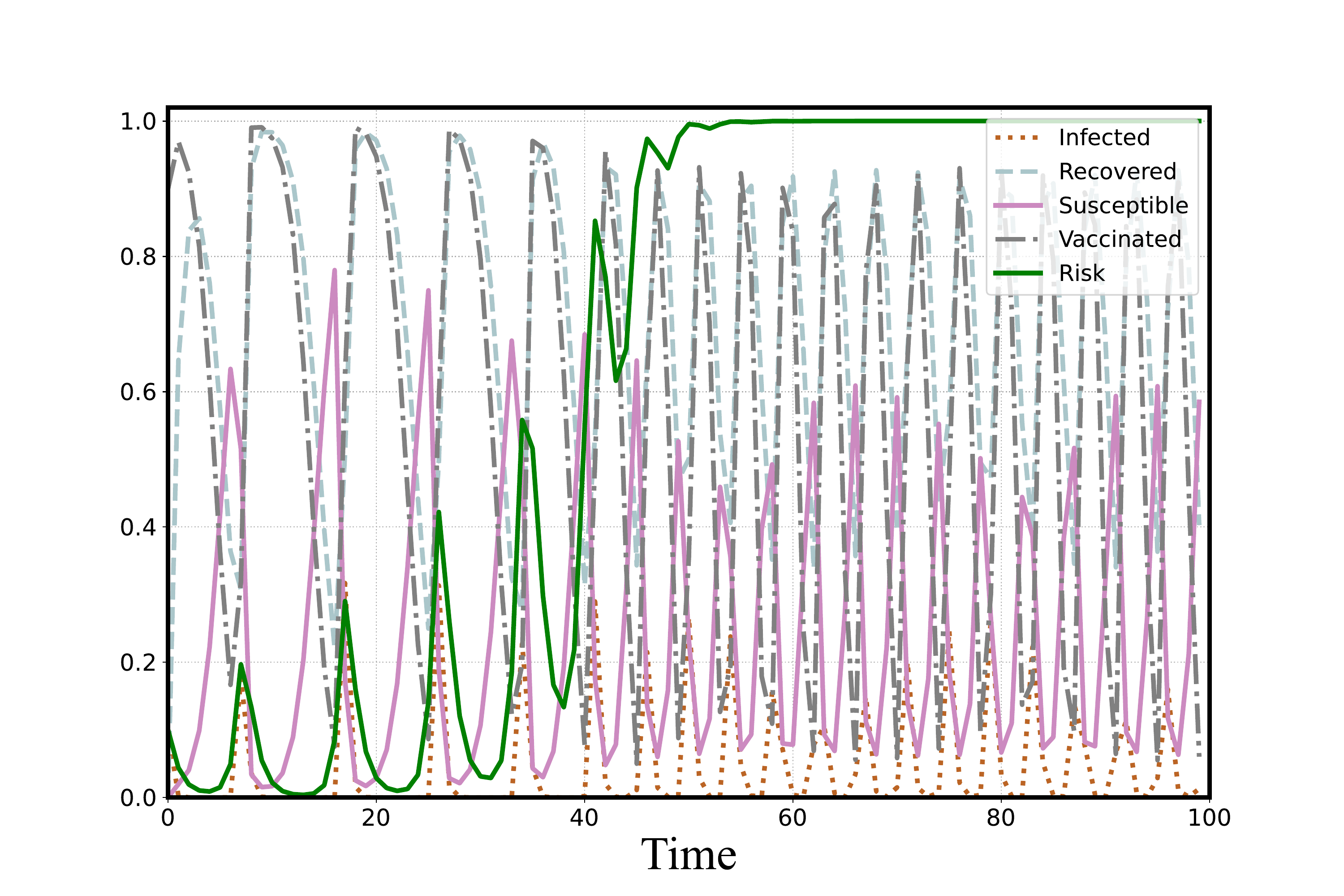}[south]{$x_0 = 0.9, n_0 = 0.1$}{SIRb-1} 
	\end{minipage}
	\begin{minipage}{0.45\linewidth} 
	\MySubFigure{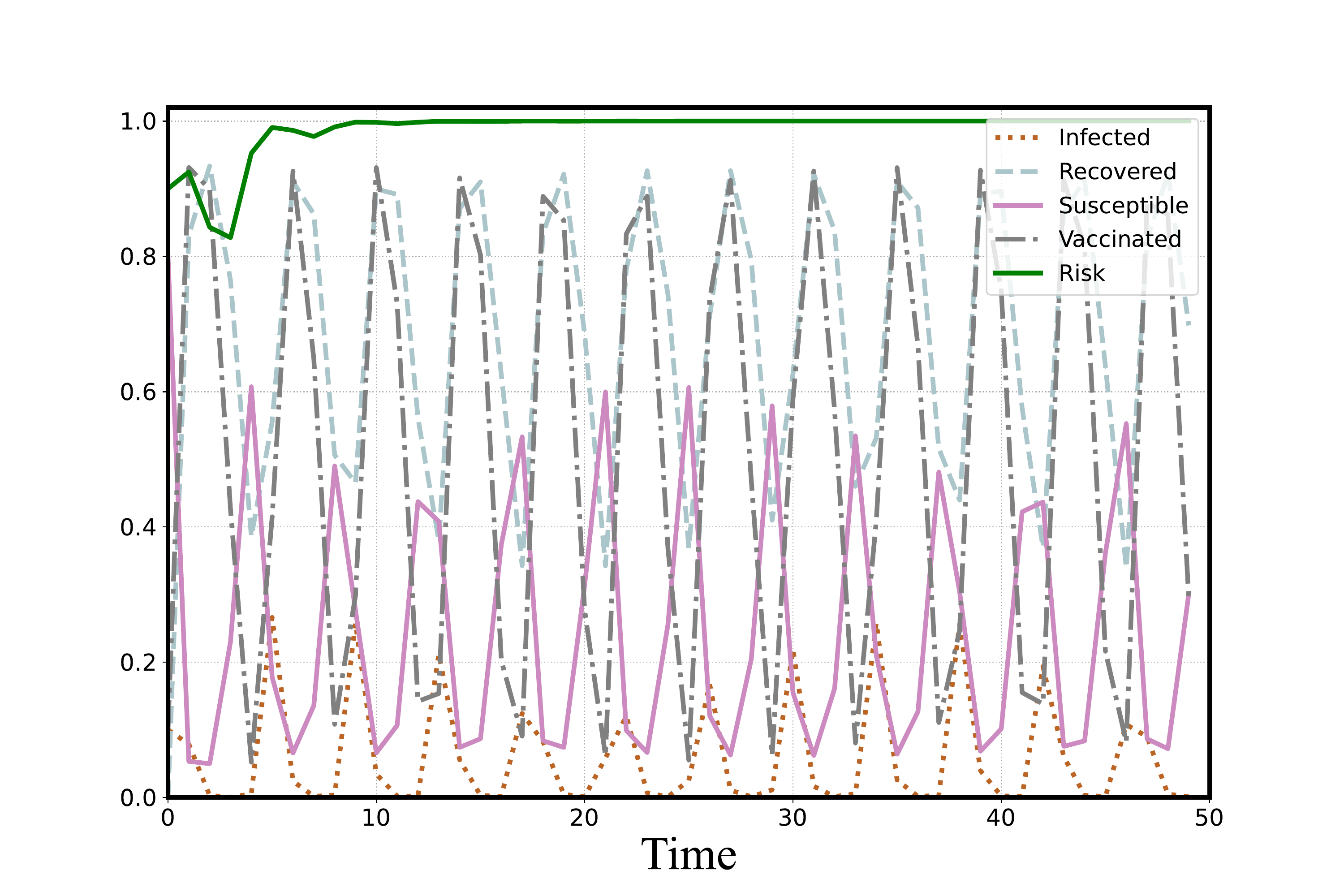}[south]{$x_0 = 0.1, n_0 = 0.9$}{SIRb-2} 
	\end{minipage}
	\addtocounter{figure}{-1} 
	\caption{\label{SIRb} Emergent vaccination scare. (a) Initial vaccination level $x_0 = 0.9$ and perceived vaccination risk $n_0=0.1$, and the final vaccination level reach a higher proportion and perceived risk reach zero. 
	(b) Initial vaccination level $x_0=0.1$ and perceived vaccination risk $n_0 = 0.9$, and the final vaccination level reach a lower proportion and perceived risk reach one. 
	For the two cases, the final perceived risk reach one but vaccination level is unstable. 
	The phenomenon of oscillating vaccination level and epidemics does not find in the model with time scale separation. 
	Parameters: $\epsilon_1 = \epsilon_2 = 0.99, \mu=1, \beta=16, \gamma=3, \theta=1, C=10, V_H=3, V_L=1$.
}
\end{figure*}
It implies that it can lead to vaccine scare if epidemic does not die out or new infected individuals increase frequently.
A similar case is COVID-19, which has been coming back and forth for more than two years, with infected individuals still present and not yet reaching the endpoint. 
Thus, it can lead to vaccine scare if the COVID-19 has been slow to die out.

The above analysis is based on the same time scale between vaccination behavior and perceived risk. 
In fact, vaccination behavior and perceived vaccination risk are not necessarily to evolve at the same time scale, that is, $\epsilon_1 \neq \epsilon_2$.
However, we find that a slight change in $\epsilon_2$ does not change the dynamical results, no matter $\epsilon_1 = 0.99$ or $\epsilon_1=0.01$.
We take $\epsilon_1=0.01$ as an example.
There are also two stable regimes when the basic reproductive ratio is moderate, regardless of $\epsilon_2=0.005$ or $\epsilon_2=0.05$ (See Fig.~\ref{SIRc}).
\begin{figure*}
	\stepcounter{figure} 
	\begin{minipage}{0.45\linewidth} 
	\MySubFigure{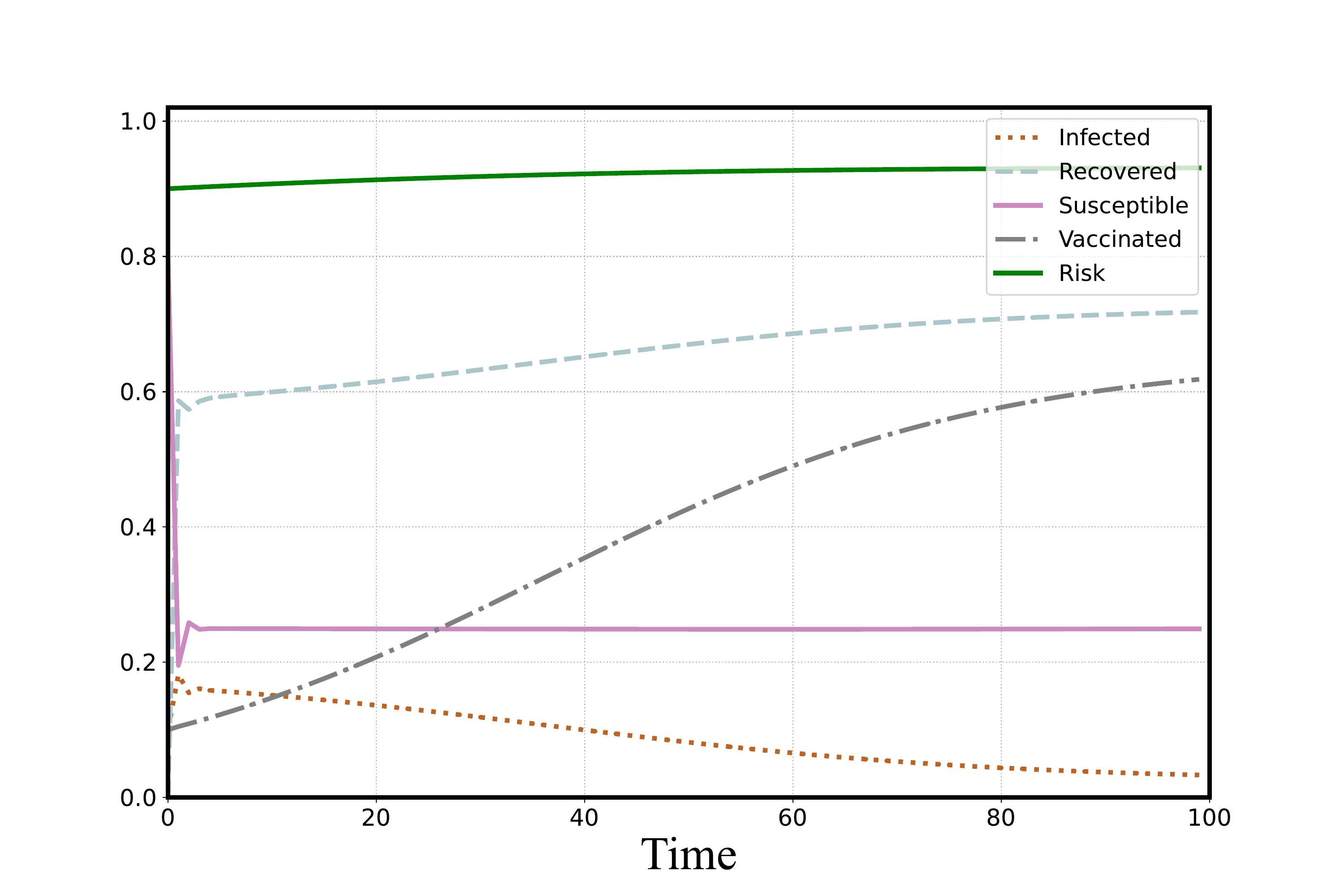}[south]{$\epsilon_1=0.01$, $\epsilon_2=0.005$}{SIRc-1} 
	\end{minipage}
	\begin{minipage}{0.45\linewidth} 
	\MySubFigure{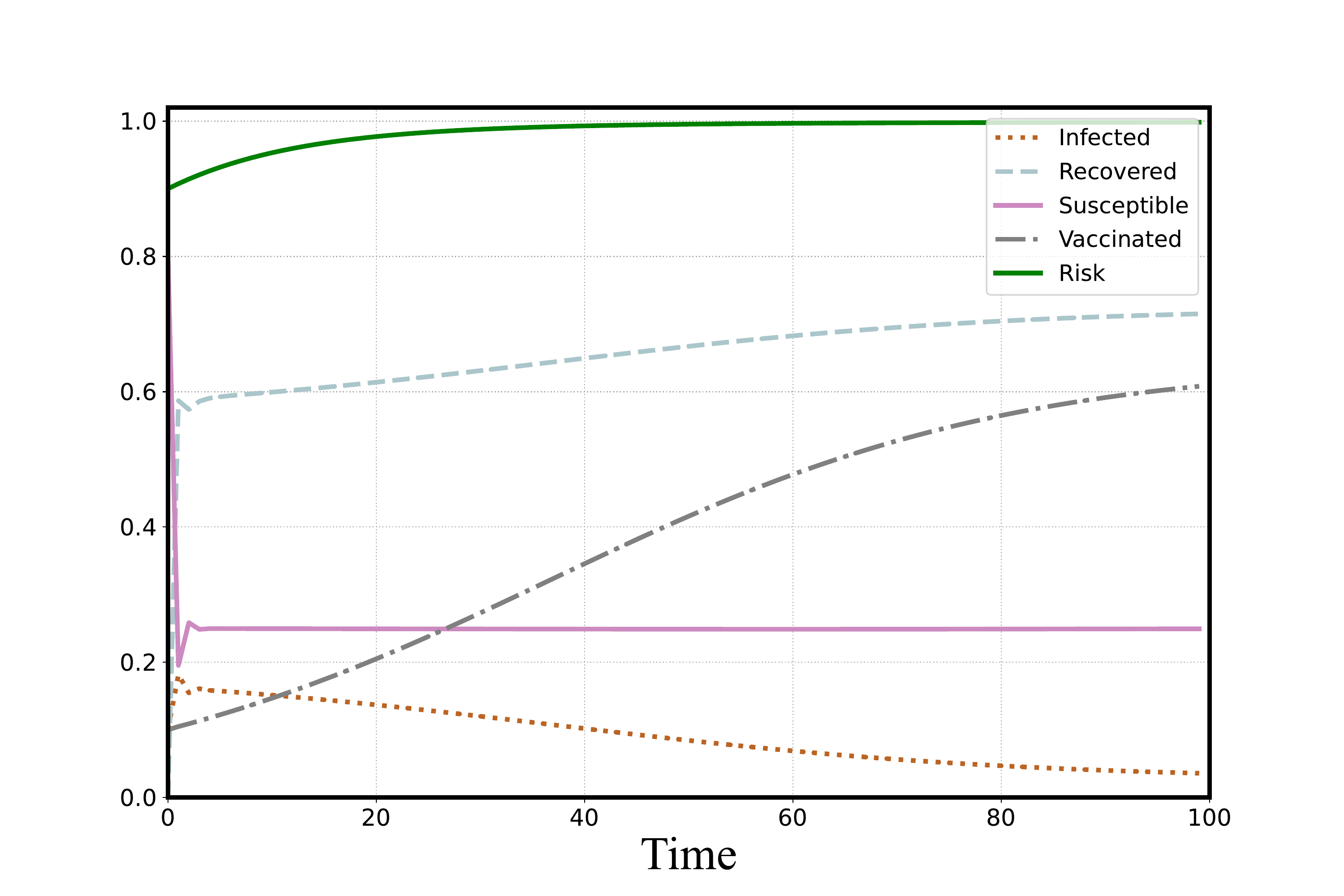}[south]{$\epsilon_1=0.01$, $\epsilon_2=0.05$}{SIRc-3} 
	\end{minipage}
	\begin{minipage}{0.45\linewidth} 
	\MySubFigure{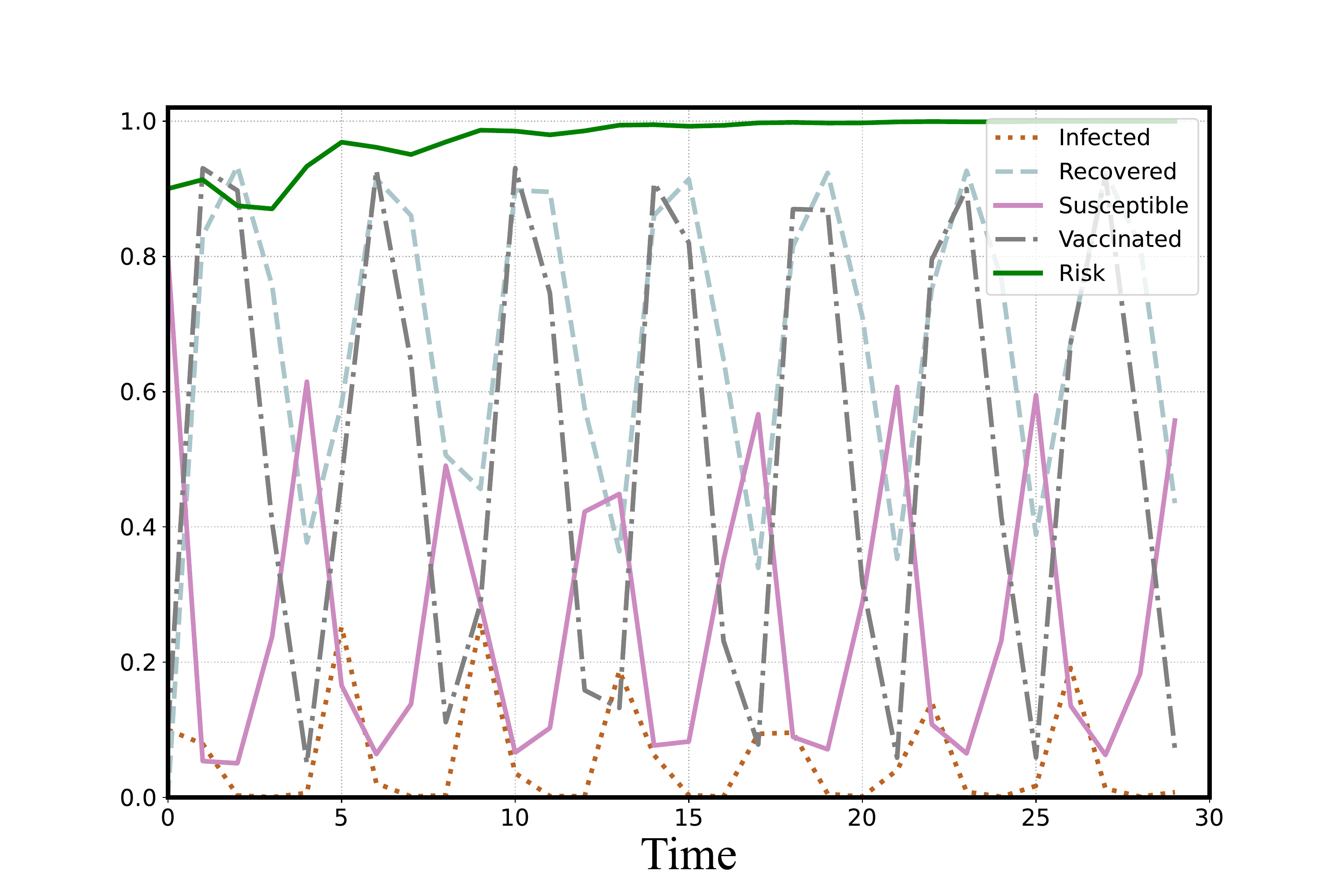}[south]{$\epsilon_1=0.99$, $\epsilon_2=0.5$}{SIRc-2} 
	\end{minipage}
	\begin{minipage}{0.45\linewidth} 
	\MySubFigure{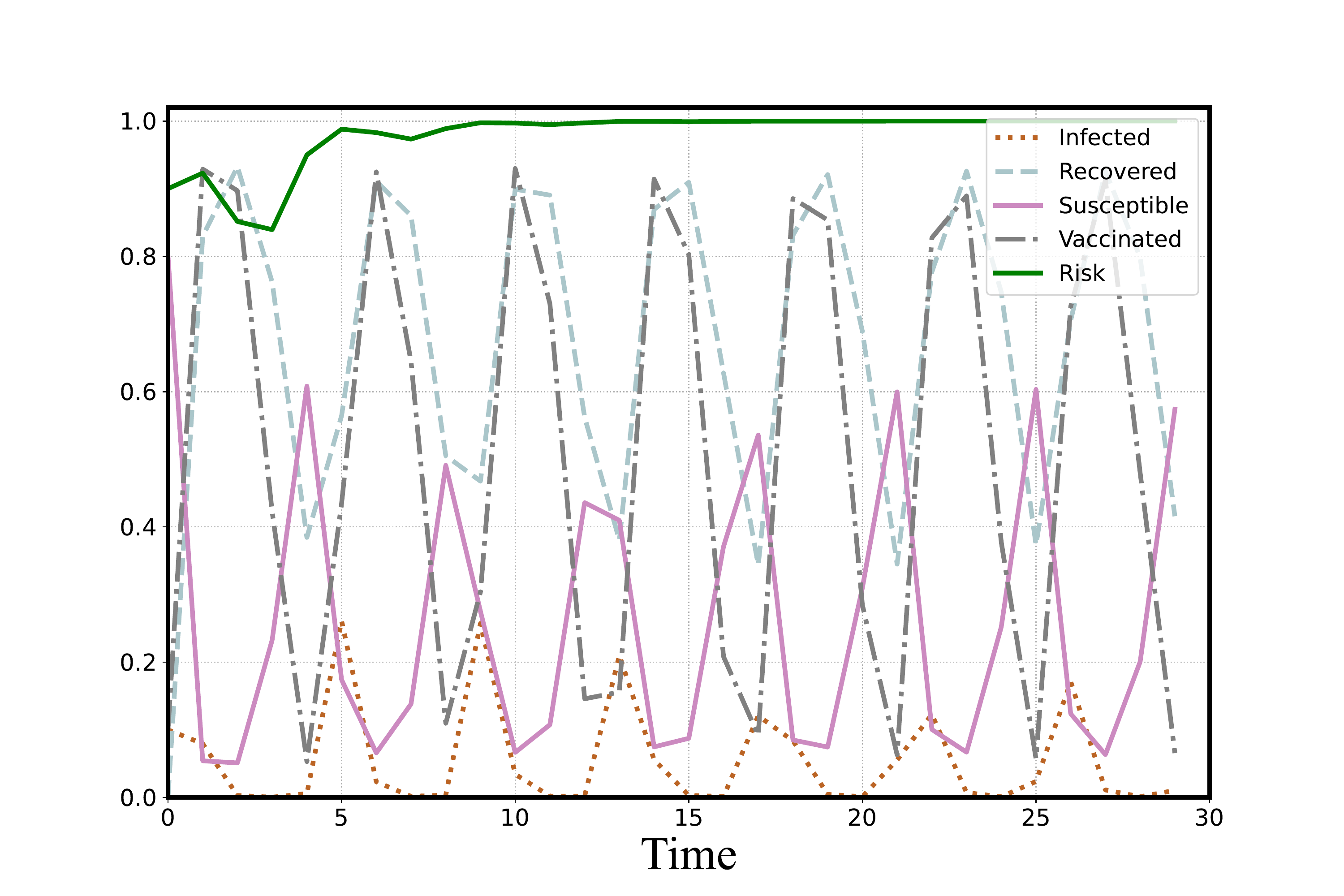}[south]{$\epsilon_1=0.99$, $\epsilon_2=0.9$}{SIRc-4} 
	\end{minipage}
	\addtocounter{figure}{-1} 
	\caption{\label{SIRc} The time scales can not change the dynamical results qualitatively. 
	When $\epsilon_1=0.01$, the final state of dynamics remains unchanged under condition (a) $\epsilon_2=0.005$ and (b) $\epsilon_2=0.05$, even there is a slight change in evolution. 
	It also remains unchanged under condition (c) $\epsilon_2=0.5$ and (d) $\epsilon_2=0.9$ when $\epsilon_1=0.99$.  
	Parameters: $\mu=1, \beta=16, \gamma=3, \theta=1, C=10, V_H=3, V_L=1$.
}
\end{figure*}

\section{Epidemic Control via Side-effect Bias}
The co-evolutionary dynamics of cognition and behavior is taken into account.
Thus, it is important to study how to promote vaccination behavior by controlling side-effect bias $\theta$, which measures the ratio of the enhancement rates to degradation rates of vaccinated and unvaccinated individuals.
We focus on a simple control rule that the side-effect bias reduces as long as the infected individuals exceed a critical value. 
Noteworthily, the side-effect bias could be adjusted by policy modification. 
For example, reporting the positive utility and data of vaccinated individuals is an effective way to decrease the side-effect bias. 
We assume that side-effect bias $\theta$ decreases if epidemic reaches a critical size, that is, $\theta=0.3$ if the proportion of infected individuals exceeds $I_t=0.001$ at any time $t$ (See Fig.~\ref{SIRd}).
\begin{figure*}
	\stepcounter{figure} 
	~~~~~~~~without control~~~~~~~~~~~~~~~~~~~~~~~~~~~~~~~~~~
	~~~~~~~~~~~~~~~~~~~~control\\
	$\epsilon = 0.1$
	\begin{minipage}{0.42\linewidth} 
	\includegraphics[width=1\linewidth]{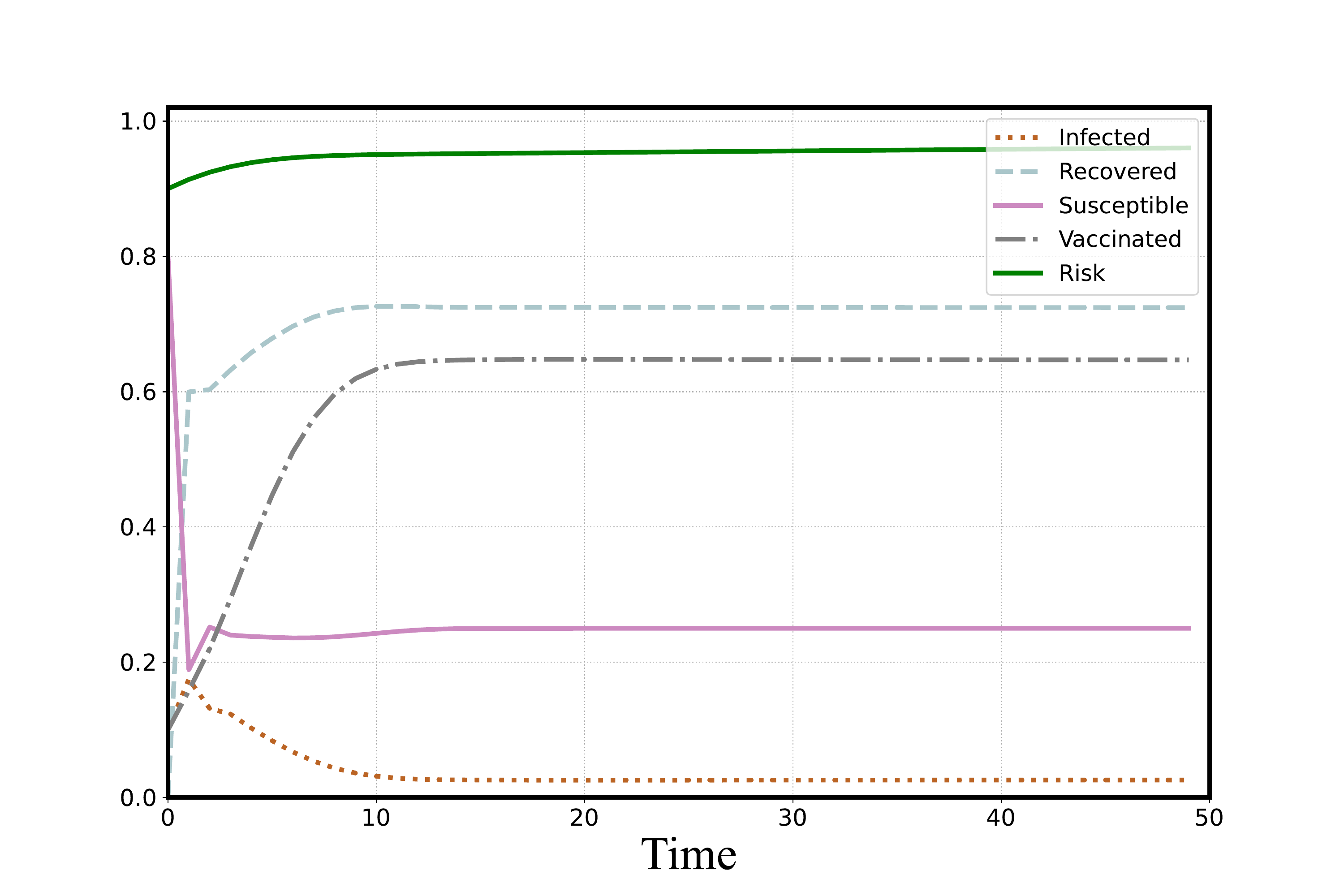}
	\end{minipage}
	\begin{minipage}{0.42\linewidth} 
	\includegraphics[width=1\linewidth]{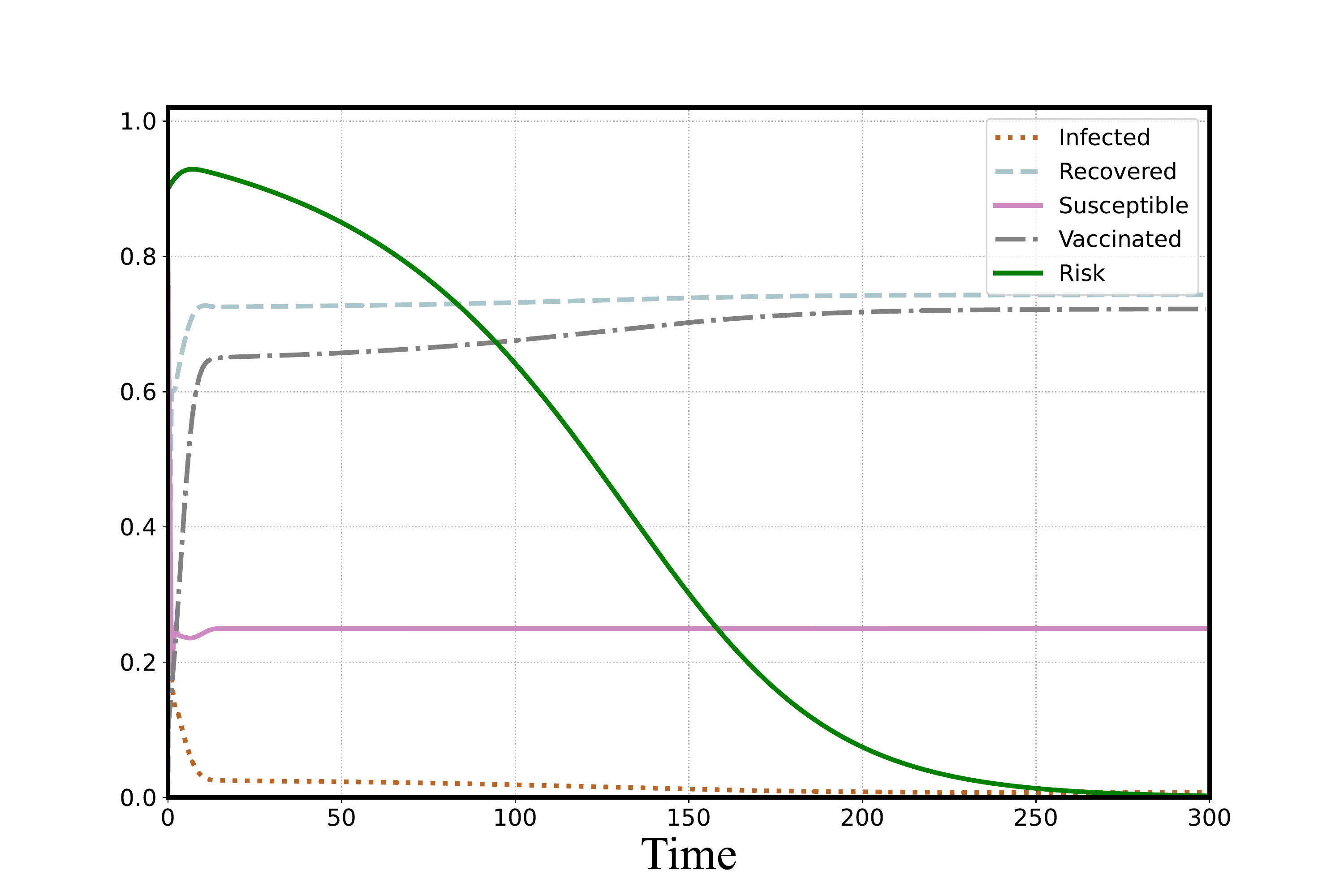}
	\end{minipage}\\
	$\epsilon = 0.5$
	\begin{minipage}{0.42\linewidth} 
	\includegraphics[width=1\linewidth]{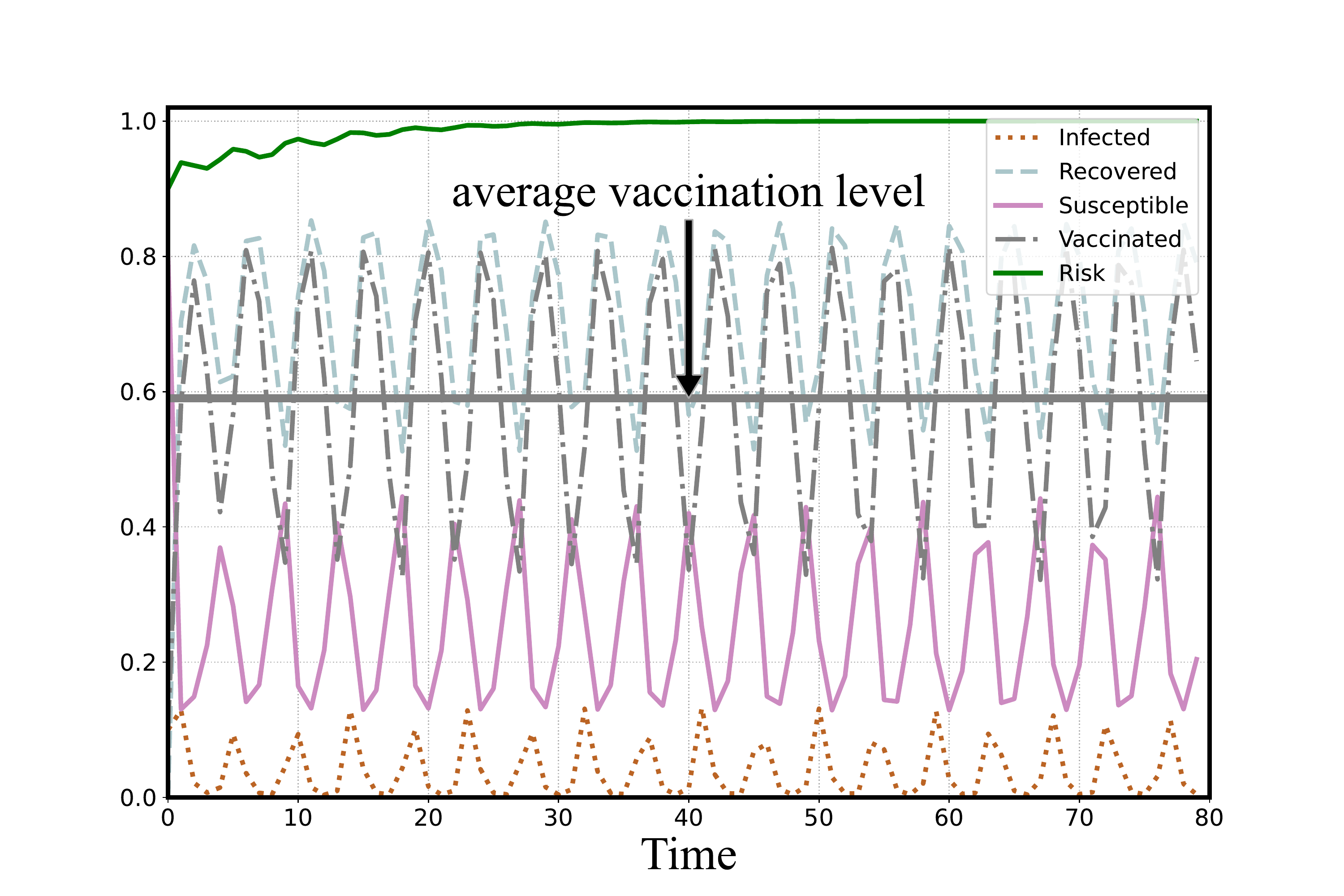}
	\end{minipage}
	\begin{minipage}{0.42\linewidth} 
	\includegraphics[width=1\linewidth]{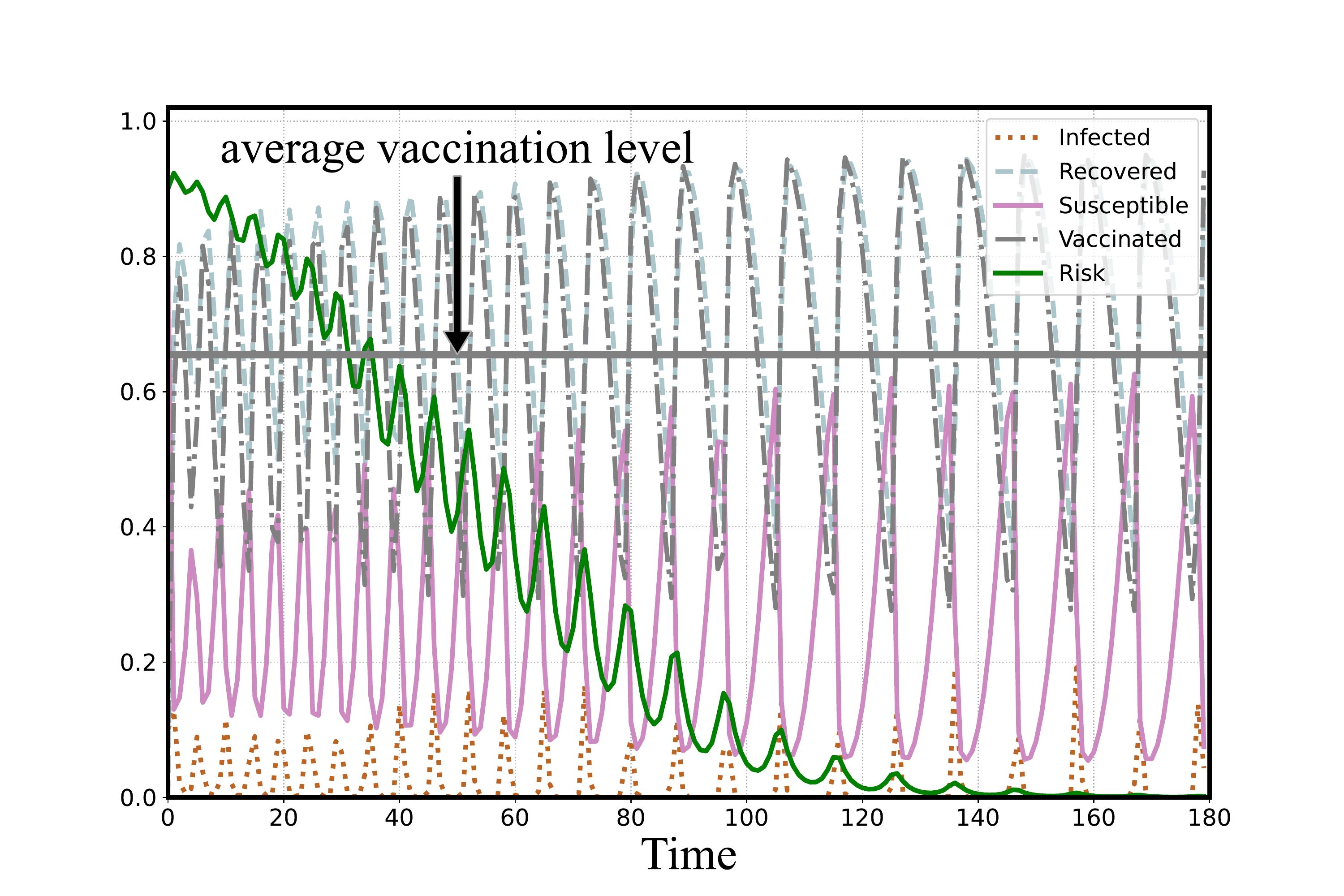}
	\end{minipage}
	\addtocounter{figure}{-1} 
	\caption{\label{SIRd} Controlling side-effect bias promotes vaccination behavior. 
	Gray dotted lines denote vaccination level, green lines denote perceived vaccination risk in four figures.
	Gray solid lines in the second row denote average vaccination level.
	The time scale of the first row is $\epsilon_1=\epsilon_2=0.1$ and that of the second row is $\epsilon_1=\epsilon_2=0.5$. 
	Controlling side-effect bias is that change $\theta$ from $1$ to $0.3$ if the proportion of infected individuals reaches $I_t=0.001$ at any time.
	The first column have no control, which implies that $\theta$ is always equal to $1$.
	And second column have controls. 
	Vaccination level increases and perceived risk decreases with controlling side-effect bias $\theta$. 
	 Parameters: $\mu=1, \beta=16, \gamma=3, \theta=1, C=10, V_H=3, V_L=1$.
}
\end{figure*}
It implies that the effect of unvaccinated on perceived vaccination risk reduces as long as the infected individuals exceed one thousand in a population with one million individuals.
It is found that it promotes the vaccination behavior at various time scales when the epidemic reaches a critical size.

Then we consider a strict condition that side-effect bias $\theta$ decreases to $0.0001$ if the proportion of infected individuals exceed $I_t=0.0001$ at any time $t$, that is, the effect of unvaccinated on perceived vaccination risk is much lower when the infected individuals exceeds one hundred. (See Fig.~\ref{SIRe}).
It is found that the vaccination behavior is also promoted and evolves much faster under this condition.
\begin{figure*}
	\stepcounter{figure} 
	\begin{minipage}{0.45\linewidth} 
	\MySubFigure{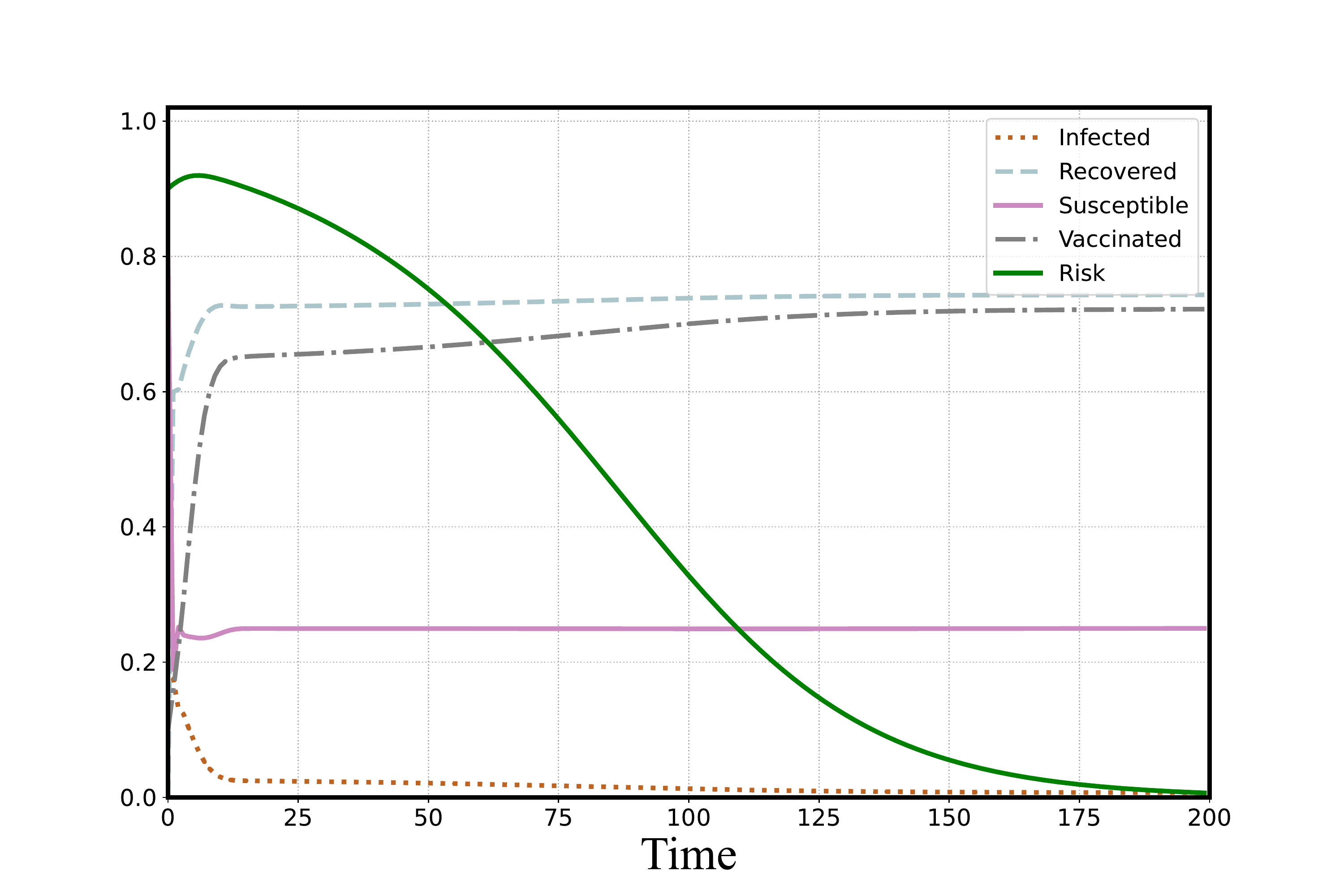}[south]{$\epsilon_1=\epsilon_2=0.1$}{SIRd-4} 
	\end{minipage}
	\begin{minipage}{0.45\linewidth} 
	\MySubFigure{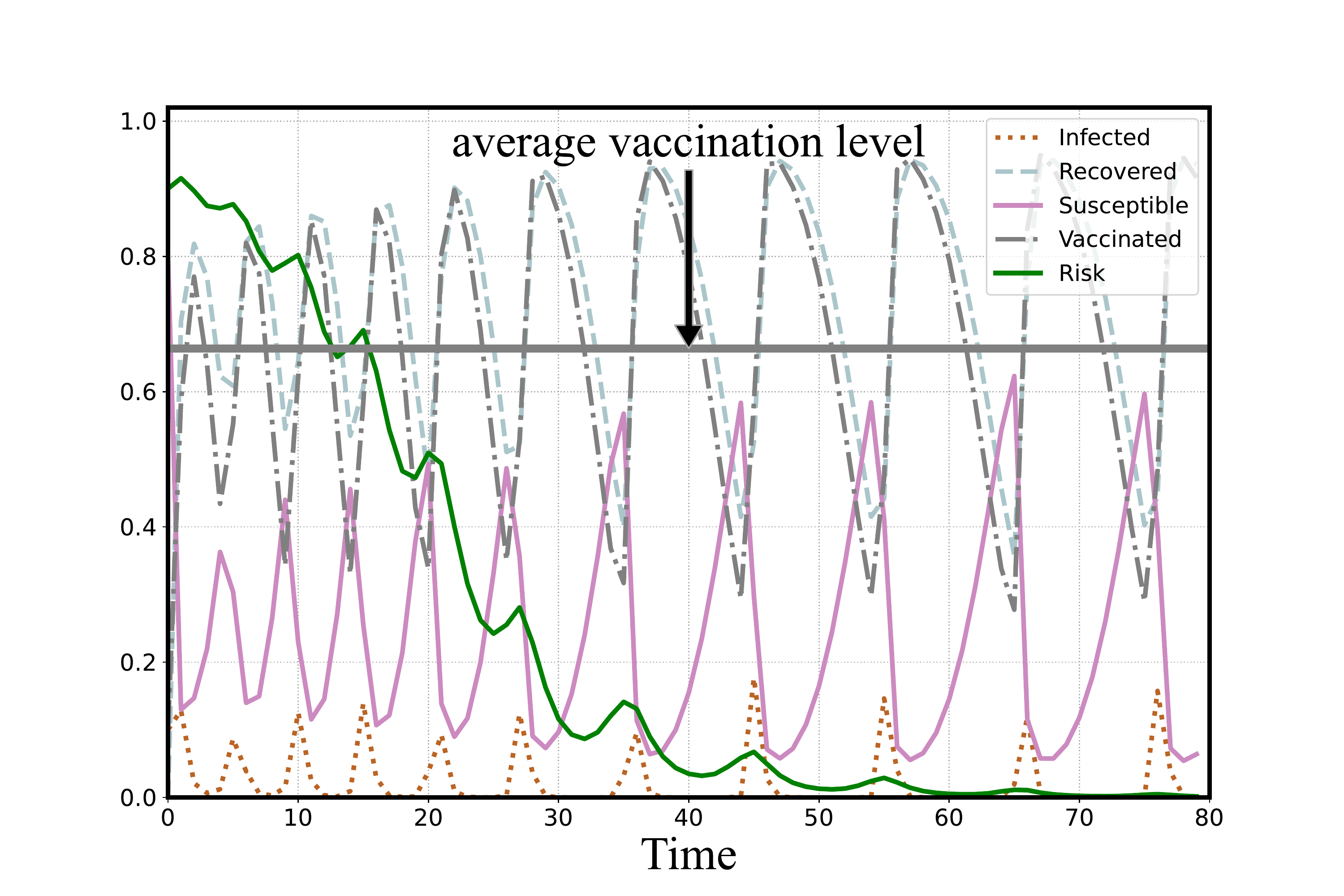}[south]{$\epsilon_1=\epsilon_2=0.5$}{SIRd-2} 
	\end{minipage}
	\addtocounter{figure}{-1} 
	\caption{\label{SIRe} Strict conditions accelerate the evolution of vaccination behavior.
	The time scale of (a) is $\epsilon_1=\epsilon_2=0.1$ and that of (b) is $\epsilon_1=\epsilon_2=0.5$. 
	Controlling side-effect bias is that change $\theta$ from $1$ to $0.0001$ if the proportion of infected individuals reach $I_t=0.0001$ at any time.
	Parameters: $\mu=1, \beta=16, \gamma=3, \theta=1, C=10, V_H=3, V_L=1$.
}
\end{figure*}

\begin{figure}
	\includegraphics[scale=0.25]{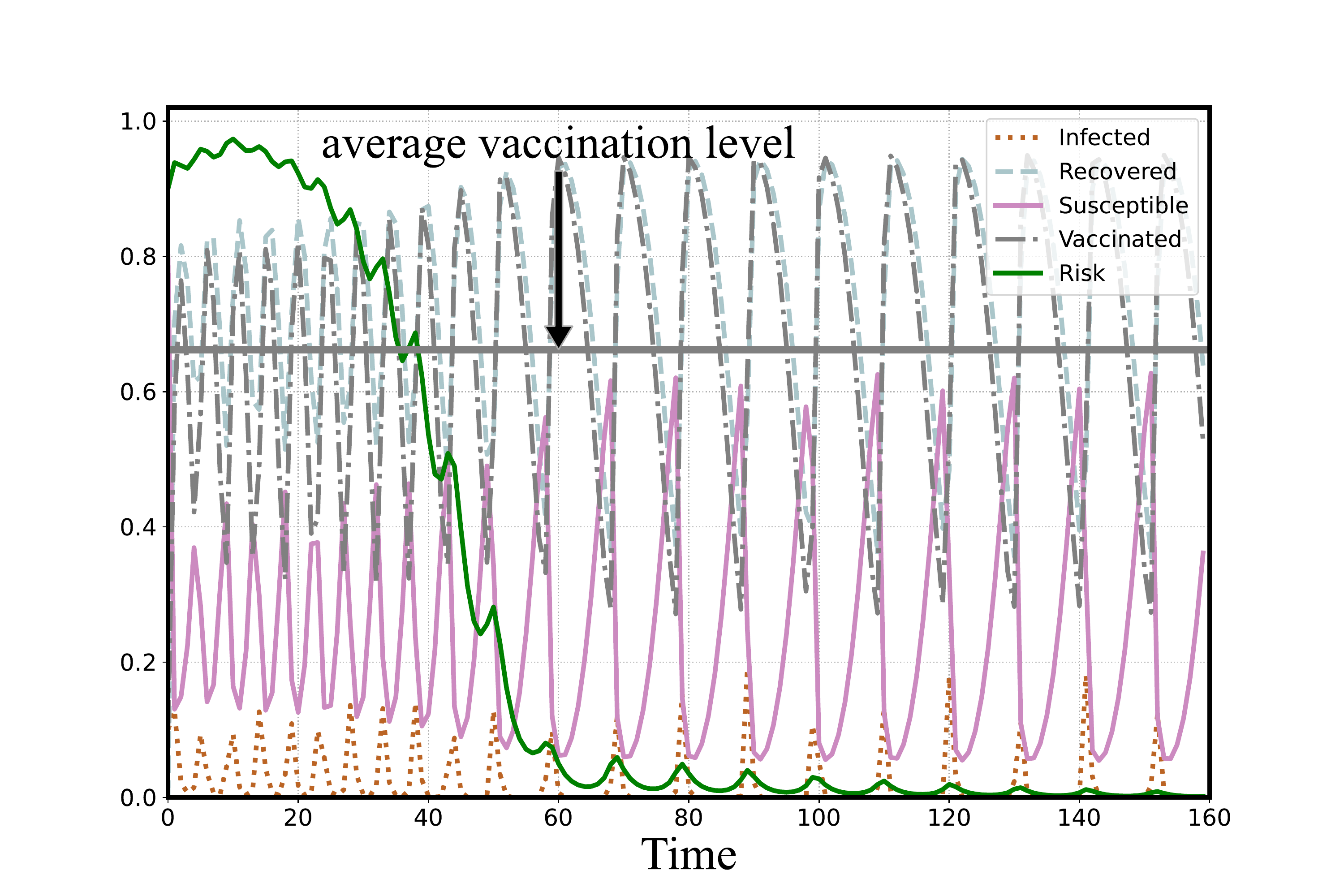}
	\caption{\label{SIRf} Vaccination behavior is promoted by controlling side-effect bias in a period of time.
	 $\theta=0.01$ if time $t \in (10, 60)$, and $\theta=1$ in other time.
	Parameters: $R_0=3.5, C=10, V_H=3, V_L=1, \theta=1,\epsilon_1=\epsilon_2=0.5$. 
}
\end{figure}

In addition, we consider one-shot controlling side-effect bias, that is, $\theta=0.01$ if time $t$ is within the interval $(10, 60)$ (See Fig.~\ref{SIRf}).
And we find that it also promotes the vaccination behavior in this case.

To sum up, we find that it is effective to promote vaccination behavior by controlling side-effect bias at both fast and slow time scales.
Thus, what we should do is to increase the positive utility reports of vaccinated individuals, and to reduce the spread of exaggerated and untrue statements of unvaccinated individuals when epidemic is present.

\section{Conclusion}
We have introduced a co-evolutionary model of epidemic, vaccination behavior and perceived vaccination risk: 
these three dynamics evolve with different time scales.
We have assumed that the increase of vaccination level inhibits the rise of vaccination risk (popular behaviors are likely to be taken as safe), and the increase of vaccination risk inhibits the rise of perceived vaccination level (high cost behavior is not likely to evolve).

We have shown how to make the whole population reach a higher vaccination level for moderate infections disease if epidemic evolves the fastest among these three dynamics:
Firstly, the lower side-effect bias is, the larger the attraction basin of the point (1) is.
The side-effect bias strengthens the role of the proportion of unvaccinated in increasing the vaccination risk.
Therefore, the side-effect bias should be as small as possible in order to make the overall vaccination level higher.
That's to say, the side-effect bias should be reduced, so as to avoid the reduction of vaccination level due to the enhanced of vaccination risk.
It suggests that, officials should adopt policies to strengthen the trust in vaccines and punish individuals who maliciously spread false information.
Secondly, the smaller the vaccination cost is the larger the attraction basin of the point is.
This is because low perceived vaccination cost promotes vaccination.
It implies that the interaction between vaccination level and vaccination risk increases the vaccination level.
Therefore, in order to make the whole population reach high vaccination level, the vaccination cost should be reduced, including time consumption, side effects etc.
Finally, the larger the basic reproductive ratio is, the larger the attraction basin of the point is.
This is because individuals recognize that the vaccination risk is obviously lower than the infection risk. 
Rational individuals tend to vaccinate, which improves the overall vaccination level.
Thus, it is effective to reduce the side-effect bias and vaccination cost when infectious diseases are not serious. 
Side-effect bias refers to the interaction between vaccination level and perceived vaccination risk.
These ways effectively prevent and control the spread of infectious diseases by promoting the vaccination and simultaneously reducing the perceived vaccination risk.
We have supposed that the promotion effect of non vaccination level on perceived risk is stronger than that of vaccination level.
However, it is challenging to measure side-effect bias.
For example, individuals who cannot vaccinate due to illness are different from those who believe and even spread false information about vaccination \cite{hajure2021}. 
The two types have different impact on the perceived vaccination risk.
It is worthy to study in the future.

Generally, we have considered these three dynamics evolve with a similar time scale.
We have found that it can lead to vaccine scare if epidemic is not serious. 
In this case, epidemic does not die out naturally and the vaccination level can not lead to herd immunity.
Outbreaks of epidemic come and go, and the proportion of susceptible, infected, recovered and vaccinated individuals continues to fluctuate. 
But the perceived vaccination risk reaches one, which refers to vaccine scare.
A similar case is COVID-19, which has been going back and forth for nearly two years, with infected individuals still present and not yet reaching the endpoint. 
Thus, it can lead to vaccine scare if the COVID-19 has been slow to die out.
In addition, we have found that it is effective to promote vaccination behavior by controlling side-effect bias at both fast and slow time scales.

\section{Discussion}
In previous studies on voluntary vaccination behavior, individuals do what others don't do \cite{fu2011,wu2011,bauch2004,bauch2005,kabir2019,peng2016,feng2018}, which is similar to snowdrift game.
Vaccinated individuals obtain immune ability, but they pay the price of time, money and side effects for taking vaccination.
Unvaccinated individuals may be infected.
But they are protected by herd immunity when enough individuals are vaccinated, and may remain healthy without paying any price.
Rational individuals do not pay the cost of vaccination, but try to be protected by herd immunity.
Voluntary vaccination cannot reach the vaccination level, which is necessary to eradicate the epidemic.
Thus it is a dilemma to vaccinate or not.
In this case, individuals would take vaccination if few other individuals do.
Otherwise many individuals take vaccination, the herd immunity is present, thus rational individuals would not take vaccination.
In contrast with previous studies, we have found that there are two stable regimes when the basic reproductive ratio is moderate, provided that epidemic evolves much faster than vaccination behavior and perceived vaccination risk.
One stable regime has high vaccination level and low perceived vaccination risk.
The other one has low vaccination level and high perceived vaccination risk.
If the initial vaccination level is high and the perceived vaccination risk is low, individuals choose to vaccinate and think the vaccination risk is low.
On the contrary, if the initial vaccination level is low and the perceived vaccination risk is high, more individuals are not likely to vaccinate and think the perceived vaccination risk is high.
Herein all the individuals do what others do.
It implies co-evolutionary vaccination is similar to a stag-hunt game \cite{ito2017,belloc2019}.
Therefore, our co-evolutionary dynamics intrinsically transform snowdrift like game to stag-hunt like game.
Intuitively, high vaccination level would lead to low perceived risk of vaccination, and low perceived risk would promote vaccination.
Thus perceived vaccination risk acts as a catalyst to yield such stag-hunt like vaccination behavior.

There are other models of evolutionary game with environmental feedback mechanism for general common tragedy \cite{weitz2016,dong2018,hilbe2018,shao2019,tilman2020,CAO2021111088,LIU2022127309}.
Weitz et al. has proposed the replicator dynamics of evolutionary game with feedback mechanism.
They linked the payoffs of classic two-strategy symmetric matrix game with environmental feedback.
And it is proved that it cannot give rise to two stable regimes \cite{weitz2016}.
This is significantly different from our model based on vaccination game and state-dependent feedback.
It suggests that the vaccination game is intrinsically complex, which is different from the matrix game.
Thus, vaccination game can not be simply regarded as a snowdrift game.

To sum up, we have pointed out that co-evolutionary dynamics between vaccination behavior and perceived risk intrinsically transform snowdrift like game to stag-hunt like game, in which perceived vaccination risk acts as a catalyst to promote vaccination behavior.
The significance of our work is 
i) to reshape the understanding of the vaccination behavior by taking into account realistic evolving vaccination risk, 
ii) to highlight the coordination-like rather than snowdrift like dynamics, which is absent in vaccination behavior in evolutionary games. 
It could open a new avenue for the game-theoretic modeling of vaccination in the future. 

\appendix
\section{\label{dynamic equations}Dynamic Equations}
At the end of the epidemic season, there are three types of individuals in the population.
Vaccinated and healthy individuals' proportion is $x$ with payoff $- V (n)$;
Unvaccinated and infected individuals' proportion is $(1-x) f (x)$ with payoff $- C$, where $f(x)$ is the probability of infection of unvaccinated individuals when the proportion of vaccination is $x$;
Unvaccinated and healthy individuals' proportion is $(1-x) (1-f (x))$ with payoff $0$.
At this time, the individual imitates the strategy of other individuals via the Fermi update rule, that is, $i$ and $j$ are randomly selected, and the probability of $i$ learning $j$' strategy is $[1+\exp{(-\beta (\pi_j-\pi_i))}]^{-1}$.
$\pi_i$ and $\pi_j$ represent the perceived payoffs of $i$ and $j$ respectively, $\beta \ge 0$ is selection intensity, which determines how much the payoff difference affects the individual's strategy.
Then the dynamics of vaccination is given by
\begin{eqnarray*}
  \dot{x}=&&x(1-x)[(1-f(x))\tanh{(\frac{\beta}{2} (-V(n)))}\\
  &&+f(x)\tanh{(\frac{\beta}{2}(-V(n)+C))}].
\end{eqnarray*}
When selection intensity $\beta$ is sufficiently small, i.e., $\beta \to 0^+$ the equation can be simplified to
\begin{eqnarray}
  \dot{x}=x(1-x)(f(x)C-V(n)).
\end{eqnarray}

Next, we consider the interaction between vaccination behavior and the perceived vaccination risk, and suppose that the evolution of the perceived vaccination risk has the same time scale with that of the proportion of vaccinators.
At the end of the epidemic season, the vaccination level of the current season will have an impact on the individual's perceived vaccination risk during the vaccination campaign.
We follow the majority rule in opinion dynamics in the sense that the more individuals are taking vaccination, the more secure the vaccination is taken as.
This is because the social groups are likely to adopt popular opinions, and the popular opinions are taken as less harmful.
As a basic toy model, we propose that
i) the risk is between $0$ and $1$;
ii) the more vaccinators there are, the less the perceived risk of the vaccination is;
iii) the more unvaccinators there are, the more the perceived risk of the vaccination is;
iv) state-dependent feedback is proposed.
At the same time, the greater the feedback intensity $\theta$ is, the stronger role the unvaccinated play in the dynamics of the perceived vaccination risk.
Thus, the dynamics of perceived vaccination risk is given by
\begin{eqnarray}
  \dot{n}=n(1-n)(-x+(1+\theta)(1-x)),
\end{eqnarray}
where $\theta$ is side-effect bias, which measures the ratio of the enhancement rates to degradation rates of vaccinated and unvaccinated individuals, respectively.
and the logistic term $n(1-n)$ ensures that the perceived vaccination risk is restrained to $[0,1]$, as it is defined within. 
In addition, the term $(-x+(1+\theta) (1-x))$ describes the side-effect bias with the assumption that the more vaccinators there are, the less the perceived risk of the vaccination is; the more unvaccinators there are, the more the perceived risk of the vaccination is.

\section{\label{fixed points}Stable Regimes}
Let $\dot{x} = 0$, and $\dot{n} = 0$ in Eq.~(\ref{eq1}), then we can obtain seven fixed points as follows:\\
(1)$(1-\frac{C}{R_0(C-V_L)},0)$, where, $ R_0 > \frac{C}{C-V_L}$\\
(2)$(1-\frac{C}{R_0(C-V_H)},1)$, where, $ R_0 > \frac{C}{C-V_H}$\\
(3)$(\frac{\theta+1}{\theta+2},-\frac{C(2-R_0+\theta)+R_0V_L}{R_0(V_H-V_L)})$, where, $ \frac {(2+ \theta) C} {C-V_L} < R_0 < \frac {(2+ \theta) C} {C-V_H}$\\
(4)(0,0)\\
(5)(0,1)\\
(6)(1,0)\\
(7)(1,1)

Since $x, n \in [0, 1]$, the first three fixed points exist only if their corresponding conditions are satisfied.
We investigate the fixed points and analyze their stabilities with the aid of Jacobian matrices: \\
i) if determinant of the matrix $det(J(x^*,n^*))$ is positive and trace of it $tr(J(x^*,n^*))$ is negative, the fixed point $(x^*,n^*)$ is stable;\\
ii)  if determinant $det(J(x^*,n^*))$ is positive and trace $tr(J(x^*,n^*))$ is positive, the fixed point $(x^*,n^*)$ is unstable;\\
iii) if determinant $det(J(x^*,n^*))$ is negative and trace $tr(J(x^*,n^*))$ is either zero or indeterminacy, the fixed point $(x^*,n^*)$ is saddle.\\
For every fixed point, the Jacobian matrix can be calculated by the following equations:
\begin{eqnarray}
J(x^*,n^*)=
\begin{pmatrix}
\frac{\partial \dot{x}} {\partial x}|_{(x^*,n^*)} & \frac{\partial \dot{x}} {\partial n}|_{(x^*,n^*)} \vspace{1ex}\\
\frac{\partial \dot{n}} {\partial x}|_{(x^*,n^*)} & \frac{\partial \dot{n}} {\partial n}|_{(x^*,n^*)}
\end{pmatrix}.
\end{eqnarray}

Here, if $x<1-\frac {1} {R_0}$,
\begin{small}
\begin{eqnarray*}
&&\frac{\partial \dot{x}} {\partial x}=C (1 - \frac {1} {R_0} - 2 x) + (n (V_H - V_L) + V_L) (-1 + 2 x)\\
&&\frac{\partial \dot{x}} {\partial n}=x(V_H - V_L)(x - 1)\\
&&\frac{\partial \dot{n}} {\partial x}=n(\theta + 2)(n - 1)\\
&&\frac{\partial \dot{n}} {\partial n}=(-1 + 2 n) (-1 - \theta + 2 x + \theta x)
\end{eqnarray*}
\end{small}
and if $x \ge 1-\frac {1} {R_0}\vspace{1ex}$,
\begin{eqnarray*}
&&\frac{\partial \dot{x}} {\partial x}=(n (V_H - V_L) + V_L) (-1 + 2 x)\\
&&\frac{\partial \dot{x}} {\partial n}=x(V_H - V_L)(x - 1)\\
&&\frac{\partial \dot{n}} {\partial x}=n(\theta + 2)(n - 1)\\
&&\frac{\partial \dot{n}} {\partial n}=(-1 + 2 n) (-1 - \theta + 2 x + \theta x)
\end{eqnarray*}

Thus, the corresponding Jacobian matrices of each fixed point can be given by following equations:\\
(1)$(1-\frac{C}{R_0(C-V_L)},0)$
\begin{equation}
J=
\begin{pmatrix}
-C (1 - \frac{1}{R_0}) + V_L & \frac{C (V_H - V_L) (C - C R_0 + R_0 V_L)}{R_0^2 (C - V_L)^2}\\
0 & \frac{C (2 + \theta - R_0) + R_0 V_L} {R_0 (C - V_L)}
\end{pmatrix}\notag
\end{equation}
The fixed point (1) exists if the condition $ R_0 > \frac{C}{C-V_L}$ is satisfied. 
When $ R_0 > \frac{(2 + \theta) C}{C-V_L}$, the fixed point (1) is stable.
When $  \frac{C}{C-V_L}< R_0 < \frac{(2 + \theta) C}{C-V_L}$, the fixed point (1) is unstable.
(2)$(1-\frac{C}{R_0(C-V_H)},1)$
\begin{equation}
J=
\begin{pmatrix}
-C (1 - \frac{1}{R_0}) + V_H & \frac{C (V_H - V_L) (C - C R_0 + R_0 V_H)}{R_0^2 (C - V_H)^2}\\
0 & -\frac{C (2 + \theta - R_0) + R_0 V_H} {R_0 (C - V_H)}
\end{pmatrix}\notag
\end{equation}
The fixed point (2) exists if the condition $ R_0 > \frac{C}{C-V_H}$ is satisfied. 
When $ \frac{C}{C-V_H}< R_0 < \frac{(2 + \theta) C}{C-V_H}$, the fixed point (2) is stable.
When $ R_0 > \frac{(2 + \theta) C}{C-V_H}$, the fixed point (2) is unstable.
(3)$(\frac{\theta+1}{\theta+2},-\frac{C(2-R_0+\theta)+R_0V_L}{R_0(V_H-V_L)})$
\begin{tiny}
\begin{eqnarray*}
J=
\begin{pmatrix}
-\frac{(1 + \theta) C}{R_0} & -\frac{(1 + \theta) (V_H - V_L)}{(2 + \theta)^2}\\
\frac{(2 + \theta) (C (2 + \theta - R_0) + R_0 V_H) (C (2 + \theta - R_0) + R_0 V_L)} {R_0^2 (V_H - V_L)^2} & 0
\end{pmatrix}\notag
\end{eqnarray*}
\end{tiny}\\
The fixed point (3) exists if the condition $ \frac {(2+ \theta) C} {C-V_L} < R_0 < \frac {(2+ \theta) C} {C-V_H}$ is satisfied. 
Since the determinant of Jacobian matrix is always negative, the fixed point (3) is saddle.
(4)(0,0)
\begin{equation}
J=
\begin{pmatrix}
C(1-\frac{1}{R_0})-V_L & 0\\
0 & 1+\theta
\end{pmatrix}\notag
\end{equation}
(5)(0,1)
\begin{equation}
J=
\begin{pmatrix}
C(1-\frac{1}{R_0})-V_H & 0\\
0 & -(1+\theta)
\end{pmatrix}\notag
\end{equation}
(6)(1,0)
\begin{equation}
J=
\begin{pmatrix}
V_L & 0\\
0 & -1
\end{pmatrix}\notag
\end{equation}
(7)(1,1)
\begin{equation}
J=
\begin{pmatrix}
V_H & 0\\
0 & 1
\end{pmatrix}\notag
\end{equation}
The last three fixed points are always exist based on our model. 
The fixed point (5) is stable if $R_0 < C/(C - V_H)$ and other points are always unstable.

After obtaining the Jacobian matrices at these fixed points, we can carry out the following stability analysis:

$\bm{Case~1:}$
\begin{eqnarray}\label{eq3}
  \frac{(2+\theta)C}{C-V_L}>\frac{C}{C-V_H},
\end{eqnarray}
which is required condition one $C\leq 2V_H-V_L$, $\theta >\frac{-C+2V_H-V_L}{C-V_H}$, or two $C>2V_H-V_L$,

1) If
\begin{equation}
  0 < R_0 < \frac{C}{C-V_L}, \notag
\end{equation}
there are four fixed points (4)(5)(6)(7), and only point (5) is stable.

2) If
\begin{equation}
  \frac{C}{C-V_L} < R_0 < \frac{C}{C-V_H}, \notag
\end{equation}
there are five fixed points (1)(4)(5)(6)(7), and only point (5) is stable.

3) If
\begin{equation}
  \frac{C}{C-V_H} < R_0 < \frac{(2+\theta)C}{C-V_L}, \notag
\end{equation}
there are six fixed points (1)(2)(4)(5)(6)(7), and only point (2) is stable.

4) If
\begin{equation}
  \frac{(2+\theta)C}{C-V_L} < R_0 < \frac{(2+\theta)C}{C-V_H}, \notag
\end{equation}
there are seven fixed points (1)(2)(3)(4)(5)(6)(7), both points (1) and (2) are stable.

5) If
\begin{equation}
  R_0 > \frac{(2+\theta)C}{C-V_H}, \notag
\end{equation}
there are six fixed points (1)(2)(4)(5)(6)(7), and only point (1) is stable.

$\bm{Case~2:}$
\begin{equation}\label{8}
  \frac{(2+\theta)C}{C-V_L} < \frac{C}{C-V_H},
\end{equation}
which is required $C < 2V_H-V_L$ and $0< \theta < \frac{-C+2V_H-V_L}{C-V_H}$,

1) If
\begin{equation}
  0 < R_0 < \frac{C}{C-V_L}, \notag
\end{equation}
there are four fixed points (4)(5)(6)(7), and only point (5) is stable.

2) If
\begin{equation}
  \frac{C}{C-V_L} < R_0 < \frac{(2+\theta) C}{C-V_L}, \notag
\end{equation}
there are five fixed points (1)(4)(5)(6)(7), and only point (5) is stable.

3) If
\begin{equation}
  \frac{(2+\theta) C}{C-V_L} < R_0 < \frac{C}{C-V_H}, \notag
\end{equation}
there are six fixed points (1)(3)(4)(5)(6)(7), both points (1) and (5) are stable.

4) If
\begin{equation}
  \frac{C}{C-V_H} < R_0 < \frac{(2+\theta)C}{C-V_H}, \notag
\end{equation}
there are seven fixed points (1)(2)(3)(4)(5)(6)(7), both points (1) and (2) are stable.

5) If
\begin{equation}
  R_0 > \frac{(2+\theta)C}{C-V_H}, \notag
\end{equation}
there are six fixed points (1)(2)(4)(5)(6)(7), and only point (1) is stable.

\section{\label{attraction basin}Attraction Basin}
Based on the above analysis, we obtain seven fixed points and find three cases of two stable fixed points at the same time, such as 4) in case 1, 3) and 4) in case 2.
These three cases show that whatever the initial state is, the final state will reach either fixed point (1) or (2), provided that $R_0$ fulfills the conditions.

The internal fixed point (3)$(\frac{\theta+1}{\theta+2},-\frac{C(2-R_0+\theta)+R_0V_L}{R_0(V_H-V_L)})$, in which $ \frac {(2+ \theta) C} {C-V_L} < R_0 < \frac {(2+ \theta) C} {C-V_H}$, is a saddle point because its corresponding matrix has two eigenvalues: one is positive and the other is negative.
The attraction basin of the two stable fixed points on the boundary are separated by a curve passing through the internal unstable fixed point (3).
In the vicinity of the fixed point (3), the curve dividing the attraction basin can be approximated as the eigenvector corresponding to the negative eigenvalue of the Jacobian matrix of this point.
We calculate two eigenvalues of its Jacobian matrix as follows:
\begin{tiny}
\begin{widetext}
\begin{eqnarray*}
&&\lambda_1=\left [- C \left(\theta + 1\right) \left(V_H \theta + 2 V_H - V_L \theta - 2 V_L\right) - \left(- \left(V_H - V_L\right) \left(\theta + 1\right) \left(\theta + 2\right) \left(4 C^{2} R_{0}^{2} - 8 C^{2} R_{0} \theta - 16 C^{2} R_{0} - C^{2} V_H \theta^{2} - 3 C^{2} V_H \theta - 2 C^{2} V_H + C^{2} V_L \theta^{2} + 3 C^{2} V_L \theta \right. \right. \right.\\
 &&\left. \left. \left. + 2 C^{2} V_L + 4 C^{2} \theta^{2} + 16 C^{2} \theta + 16 C^{2} - 4 C R_{0}^{2} V_H - 4 C R_{0}^{2} V_L + 4 C R_{0} V_H \theta + 8 C R_{0} V_H + 4 C R_{0} V_L \theta + 8 C R_{0} V_L + 4 R_{0}^{2} V_H V_L\right) \right)^{\frac{1}{2}}\right]/
 \left[2 R_{0} \left(V_H \theta + 2 V_H - V_L \theta - 2 V_L\right) \right],\\
&&\lambda_2   = \left[- C \left(\theta + 1\right) \left(V_H \theta + 2 V_H - V_L \theta - 2 V_L\right) + \left(- \left(V_H - V_L\right) \left(\theta + 1\right) \left(\theta + 2\right) \left(4 C^{2} R_{0}^{2} - 8 C^{2} R_{0} \theta - 16 C^{2} R_{0} - C^{2} V_H \theta^{2} - 3 C^{2} V_H \theta - 2 C^{2} V_H + C^{2} V_L \theta^{2} + 3 C^{2} V_L \theta \right. \right. \right.\\
&&\left. \left. \left. + 2 C^{2} V_L + 4 C^{2} \theta^{2} + 16 C^{2} \theta + 16 C^{2} - 4 C R_{0}^{2} V_H - 4 C R_{0}^{2} V_L + 4 C R_{0} V_H \theta + 8 C R_{0} V_H + 4 C R_{0} V_L \theta + 8 C R_{0} V_L + 4 R_{0}^{2} V_H V_L\right)\right)^{\frac{1}{2}}\right]/
\left[2 R_{0} \left(V_H \theta + 2 V_H - V_L \theta - 2 V_L\right)\right],
\end{eqnarray*}
\end{widetext}
\end{tiny}
where, $\lambda_1>0$ and $\lambda_2<0$.
Then the eigenvector corresponding to negative eigenvalue $\lambda_2$ is $\boldsymbol{\eta}=[\eta_1,1]^T$, where,
\begin{tiny}
\begin{widetext}
\begin{eqnarray*}
 &&\eta_1 = - R_{0} \left(C \left(\theta + 1\right) \left(V_H \theta + 2 V_H - V_L \theta - 2 V_L\right) - \left(- \left(V_H - V_L\right) \left(\theta + 1\right) \left(\theta + 2\right) \left(4 C^{2} R_{0}^{2} - 8 C^{2} R_{0} \theta - 16 C^{2} R_{0} - C^{2} V_H \theta^{2} - 3 C^{2} V_H \theta - 2 C^{2} V_H + C^{2} V_L \theta^{2}  \right. \right. \right. \\
 && \left. \left. \left.+ 3 C^{2} V_L \theta + 2 C^{2} V_L + 4 C^{2} \theta^{2} + 16 C^{2} \theta + 16 C^{2} - 4 C R_{0}^{2} V_H - 4 C R_{0}^{2} V_L + 4 C R_{0} V_H \theta + 8 C R_{0} V_H + 4 C R_{0} V_L \theta + 8 C R_{0} V_L + 4 R_{0}^{2} V_H V_L\right)\right)^{\frac{1}{2}}\right) \left(V_H^{2} - 2 V_H V_L + V_L^{2}\right)/\\
&& \left[2 \left(V_H \theta + 2 V_H - V_L \theta - 2 V_L\right) \left(C^{2} R_{0}^{2} \theta + 2 C^{2} R_{0}^{2} - 2 C^{2} R_{0} \theta^{2} - 8 C^{2} R_{0} \theta - 8 C^{2} R_{0} + C^{2} \theta^{3} + 6 C^{2} \theta^{2} + 12 C^{2} \theta + 8 C^{2} - C R_{0}^{2} V_H \theta - 2 C R_{0}^{2} V_H - C R_{0}^{2} V_L \theta - 2 C R_{0}^{2} V_L \right. \right. \\
&&\left. \left.+ C R_{0} V_H \theta^{2} + 4 C R_{0} V_H \theta + 4 C R_{0} V_H + C R_{0} V_L \theta^{2} + 4 C R_{0} V_L \theta + 4 C R_{0} V_L + R_{0}^{2} V_H V_L \theta + 2 R_{0}^{2} V_H V_L\right)\right].
\end{eqnarray*}
\end{widetext}
\end{tiny}

In the vicinity of the fixed point (3), the curve of 4) in case 1 can be described by $f(x, n) = f(x^*, n^*) + (x^*, n^*) \boldsymbol{\eta}$ approximately.
Our purpose is to get which stable fixed point has a larger attraction basin.
Thus, for the linear approximation $n=k x+b$ for a large parameter region, i.e., $x, n \in [0,1]$, it is acceptable to evaluate which stable regime has a larger attraction basin.
The slope $k$ and intercept $b$ of the linear equation are given by:
\begin{tiny}
\begin{widetext}
\begin{eqnarray*}
  &&k=
  -\left(2 (2 + \theta)^2 (C (2 + \theta - R_0) + R_0 V_H) (C (2 + \theta - R_0) + R_0 V_L) \right)/
  \left(R_0 (V_H - V_L) \bigg((2 + 3 \theta + \theta^2) C (V_H - V_L) +
  \sqrt{(2 + 3 \theta + \theta^2) (V_H - V_L)}   \right.\\
  && \left.\times\sqrt{
   -4 R_0^2 V_H V_L - 4  C (2 + \theta -R_0) R_0 (V_H + V_L)
  + C^2 (\theta (-16 + 8 R_0 + 3 V_H - 3 V_L)
  + \theta^2 (-4 + V_H - V_L)
  - 2 (8 - 8 R_0 + 2 R_0^2 - V_H + V_L))} \bigg)\right),\\
 && b=
\left((C (2 + \theta - R_0) + R_0 V_L) \bigg(4 R_0 V_H + 6 \theta R_0 V_H + 2 \theta^2 R_0 V_H + (2 + 3 \theta + \theta^2) C (4 + 2 \theta - 2 R_0 - V_H + V_L) - \right.\\
&&\left. \sqrt{
-(2 + 3 \theta + \theta^2) (V_H - V_L) \left(4 R_0^2 V_H V_L + 4 C (2 + \theta - R_0) R_0 (V_H + V_L) + C^2 \left(\theta^2 (4 - V_H + V_L)
+ 2 (8 - 8 R_0 + 2 R_0^2 - V_H + V_L) + \theta (16 - 8 R_0 - 3 V_H + 3 V_L)\right)\right)
}\bigg)\right)/\\
&&\left(R_0 (V_H - V_L) \bigg((2 + 3 \theta + \theta^2) C (V_H - V_L) +
  \sqrt{(2 + 3 \theta + \theta^2) (V_H - V_L)}   \right.\\
  &&\left. \times\sqrt{
   -4 R_0^2 V_H V_L - 4  C (2 + \theta -R_0) R_0 (V_H + V_L)
  + C^2 \left(\theta (-16 + 8 R_0 + 3 V_H - 3 V_L)
  + \theta^2 (-4 + V_H - V_L)
  - 2 (8 - 8 R_0 + 2 R_0^2 - V_H + V_L)\right)} \bigg)\right).
\end{eqnarray*}
\end{widetext}
\end{tiny}

The existence of saddle point (3) implies that the basic reproductive number $R_0$ is not too large but guarantee that the disease spread in the individuals and do not disappear naturally, that is, $1< \frac{(2+\theta)C}{C-V_L} < R_0 < \frac{(2+\theta)C}{C-V_H}$, there are two stable fixed points in this dynamics model.
At this time, we adjust these five parameters $C$, $V_H$, $V_L$, $\theta$ and $R_0$ to stabilize as many initial points as possible at the stable fixed point (1)$(1-\frac{C}{R_0(C-V_L)},0)$ that we try to find the mechanisms to give rise to the higher vaccination level and the lower perceived vaccination risk.

\nocite{*}

\bibliography{revised}

\end{document}